\input amstex
\documentstyle{amsppt}
\document
\topmatter
\title
On the kernel and the image of the rigid analytic regulator
\endtitle
\title
On the kernel and the image of the rigid analytic regulator
in positive characteristic
\endtitle
\author
Ambrus P\'al
\endauthor
\date
May 19, 2009.
\enddate
\address Department of Mathematics, 180 Queen's Gate, Imperial College,
London SW7 2AZ, United Kingdom\endaddress
\email a.pal\@imperial.ac.uk\endemail
\footnote" "{\it 2000 Mathematics Subject Classification. \rm 19F27
(primary), 19D45 (secondary).}
\abstract
We will formulate and prove a certain reciprocity law relating certain residues
of the differential symbol $dlog^2$ from the $K_2$ of a Mumford curve to the rigid analytic regulator constructed by the author in a previous paper. We will use this result to deduce some  
consequences on the kernel and image of the rigid analytic regulator
analogous to some old conjectures of Beilinson and Bloch on the complex
analytic regulator. We also relate our construction to the symbol defined
by Contou-Carr\`ere and to Kato's residue homomorphism, and we show that
Weil's reciprocity law directly implies the reciprocity law of Anderson
and Romo. 
\endabstract
\endtopmatter

\heading 1. Introduction and announcement of results
\endheading

\definition{Motivation 1.1} In the paper [17] the author introduced
the concept of the rigid analytic regulator. This is a homomorphism
from the motivic cohomology group $H^2_{\Cal M}(X,\Bbb Z(2))$ of a Mumford curve $X$ over a local field $F$ into the $F^*$-valued harmonic cochains of the graph of components of the special fiber. It is defined through
non-archimedean integration, hence it is elementary in nature and it is
amenable to computation. In particular the author was able to compute its
value on some explicit elements of the $K_2$ of Drinfeld modular curves
constructed using modular units and relate it to special values of
$L$-functions in the paper [18]. It is quite reasonable to consider this
result as a function field analogue of Beilinson's classical theorem on
the $K_2$ of elliptic modular curves as well as the rigid analytic
regulator is a non-archimedean analogue of the complex analytic
Beilinson-Bloch-Deligne regulator.

An old conjecture of Bloch and Beilinson predicts that the complex
analytic regulator is injective on the motivic cohomology group $H^2_{\Cal M}(\goth X,\Bbb Q(2))$ of a regular integral model $\goth X$ of a smooth,
projective curve $X$ defined over a number field and the image of $H^2_{\Cal M}(\goth X,\Bbb Z(2))$ is a $\Bbb Z$-lattice of a conjecturally described rank. Hence it is desirable to understand the basic properties of the rigid
analytic regulator such as its image and kernel, partially because any
analogous result would be evidence towards the conjecture mentioned
above. We offer the following result: let $X$ be a Mumford curve which is the general fiber of a regular quasi-projective surface $\goth X$ fibred over a smooth affine curve defined over a finite field.  We show that the kernel of the rigid analytic regulator in $H^0(\goth X,K_{2,\goth X})$ is a $p$-divisible group and the image of $H^0(\goth X,K_{2,\goth X})$ is a finitely generated $\Bbb Z$-module whose closure in the $p$-completion of the target group is a $\Bbb Z_p$-module of the same rank which is at most as large as conjectured by Beilinson. In particular the kernel of this map is torsion if Parshin's conjecture holds. As a key ingredient of the proof we compute certain residues of the logarithmic differentials of elements of $H^2_{\Cal M}(\goth X,\Bbb Z(2))$ in terms of the rigid analytic regulator, generalising a formula of Osipov, proved for fields of zero characteristic in [16], to any characteristic. (A closely related result has been obtained by M. Asakura in a recent preprint [1] for certain two dimensional local fields of zero characteristic using similar methods.) Of course the case
of positive characteristic is the most intricate, due to the lack of logarithm.
This formula can be considered as a relative of the explicit reciprocity law of Kato and Besser's theorem expressing the Coleman-de Shalit regulator in terms of the syntomic regulator (see [19] and [4]), although it is simpler to prove. Then our main theorems follow from the Bloch-Gabber-Kato theorem, Deligne's purity theorem and the degeneration of the slope spectral sequence. As an important intermediate step we also relate the rigid analytic regulator to Kato's residue homomorphism for higher local
fields and the Contou-Carr\`ere symbol. The symbols mentioned above are essentially just different formalizations of the same phenomenon which was discovered
independently at least three times. As an easy application of our results
we show that the latter is bilinear and satisfies the Steinberg relation.
At the same stroke we show that Weil's reciprocity law directly implies
the reciprocity law of Anderson and Romo. 
\enddefinition
\definition{Notation 1.2} By slightly extending the usual terminology we will call a scheme $C$ defined over a field a curve if it is reduced, locally of finite type and of dimension one. A curve $C$ is said to have normal crossings if every singular point of $C$ is an ordinary double point in the usual sense. We say that a curve $C$ over a field $\bold f$ is totally degenerate if it has normal crossings, every ordinary double point is defined over $\bold f$ and its irreducible components are projective lines over $\bold f$. For any curve $C$ with normal crossings let $\widetilde C$ denote its normalization, and let $\widetilde S(C)$ denote the pre-image of the set $S(C)$ of singular points of $C$. We denote by $\Gamma(C)$ the oriented graph whose set of vertices is the set of irreducible components of $\widetilde C$, and its set of edges
is the set $\widetilde S(C)$ such that the initial vertex of any edge
$x\in\widetilde S(C)$ is the irreducible component of $\widetilde C$
which contains $x$ and the terminal vertex of $x$ is the irreducible
component which contains the other element of $\widetilde S(C)$ which
maps to the same singular point with respect to the normalization map as
$x$. The normalization map identifies the irreducible components of $C$ and $\widetilde C$ which we will use without further notice.
\enddefinition
\definition{Definition 1.3} For any (oriented) graph $G$ let $\Cal V(G)$ and $\Cal E(G)$ denote its set of vertices and edges, respectively. Let $G$ be a locally
finite oriented graph which is equipped with an involution $\overline{\cdot}:\Cal E(G)\rightarrow\Cal E(G)$ such that for each edge $e\in\Cal E(G)$ the original and terminal vertices of the edge $\overline e\in\Cal E (G)$ are the terminal and original vertices of $e$, respectively. The edge $\overline e$ is called the edge $e$ with reversed orientation. Let $R$ be a commutative group. A function $\phi:\Cal E(G)\rightarrow R$ is called a harmonic $R$-valued cochain, if it satisfies the following conditions:

\noindent $(i)$ We have:
$$\phi(e)+\phi(\overline e)=0\text{\ }(\forall e\in\Cal E(G)).$$

\noindent $(ii)$ If for an edge $e$ we introduce the notation $o(e)$ and
$t(e)$ for its original and terminal vertex respectively,
$$\sum_{{e\in\Cal E(G)\atop o(e)=v}}\phi(e)=0\text{\ }
(\forall v\in\Cal V(G)).$$
We denote by $\Cal H(G,R)$ the group of $R$-valued harmonic
cochains on $G$. 
\enddefinition
\definition{Notation 1.4} Let $\bold k$ be a perfect field and let $B$ be a smooth irreducible projective curve over $\bold k$. Let $\infty$ be a closed point of $B$ and let $F$ denote the completion of the function field of $B$ at $\infty$. Let $\Cal O$ denote the valuation ring of $F$ and let $\pi:\goth X\rightarrow B$ be a regular irreducible projective surface fibred over $B$ such that the fiber $\goth X_{\infty}$ of $\goth X$ over $\infty$ is totally degenerate. Then the base change $X$ of $\goth X$ to $F$ is a Mumford curve over $F$ which has a regular, semistable model over the spectrum of $\Cal O$ whose special fiber is $\goth X_{\infty}$. The rigid
analytic regulator introduced in [17] is a homomorphism:
$$\{\cdot\}:H^2_{\Cal M}(X,\Bbb Z(2))\rightarrow\Cal H(\Gamma(\goth X_{\infty}),F^*).$$
Let $\goth U\subset\goth X$ be an open subvariety such that its complement is a normal crossings divisor $D$ which is the preimage of a finite set of
closed points of $B$ containing $\infty$. Let the symbol $\{\cdot\}$ also
denote the composition of the functorial map $H^0(\goth U,K_{2,\goth U})
\rightarrow H^2_{\Cal M}(X,\Bbb Z(2))$ and the rigid analytic regulator.
\enddefinition
The main result of this paper is the following:
\proclaim{Theorem 1.5} Assume that $\bold k$ has characteristic $p$. Then the map:
$$\{\cdot\}:H^0(\goth U,K_{2,\goth U})/p^nH^0(\goth U,K_{2,\goth U})
\rightarrow\Cal H(\Gamma(\goth X_{\infty}),F^*/(F^*)^{p^n})$$
induced by the regulator $\{\cdot\}$ is injective for every natural number $n$.
\endproclaim
\definition{Notation 1.6} We call a $\Bbb Z$-submodule $\Lambda$ of a Hausdorff topological group $G$ which is the direct sum of a discrete group and a pro-$p$ group $p$-saturated, if it is finitely generated and the map $\Lambda\otimes\Bbb Z_p\rightarrow\widehat G$ is an injection, where
$\widehat G$ is the $p$-completion of $G$. Note that every discrete,
finitely generated $\Bbb Z$-submodule is $p$-saturated. Assume now that $\bold k$ is a finite field of characteristic $p$. The result above, Deligne's purity theorem and the degeneration of the slope spectral sequence imply the following:
\enddefinition
\proclaim{Corollary 1.7} The image of the regulator $\{\cdot\}:H^0(\goth U,K_{2,\goth U})
\rightarrow\Cal H(\Gamma(\goth X_{\infty}),F^*)$ is $p$-saturated and its rank is at most as large as the rank of the group $\Cal H(\Gamma(D),\Bbb
Z)$. This lattice is discrete if $D=\goth X_{\infty}$. The kernel of this regulator is a $p$-divisible group. In particular it is torsion if Parshin's conjecture holds for $\goth X$, and it is finite if the Bass conjecture holds for $\goth X$.
\endproclaim
\definition{Further Results 1.8} Let $A$ be a local Artinian ring with residue field $\bold k$ which is again allowed to be an arbitrary perfect field. In [7] a map:
$$\langle\cdot,\cdot\rangle:A((t))^*\times A((t))^*\rightarrow A^*$$
was defined, called the Contou-Carr\`ere symbol in [2], where $t$ is a
variable. The Contou-Carr\`ere symbol is equal to the tame symbol if $A$
is a field. In [2] G$\text{.}$ An\-der\-son and F$\text{.}$ P$\text{.}$
Romo proved that the Contou-Carr\`ere symbol is bilinear and proved a
reciprocity law for it. Their proof is the generalization of the proof of
the residue theorem by Tate and the Weil reciprocity law by Arabello, De
Concini and Kac. They work directly over Artinian rings, so they have to
develop an elaborate theory generalizing all concepts appearing in the
proofs quoted above for Artinian rings. In this article we will present a
different proof of the bilinearity of this symbol and their reciprocity
law. It is based on the observation that if $A$ is the quotient of a
discrete valuation ring, then the Contou-Carr\`ere symbol $\langle
f,g\rangle$ of any $f$ and $g$ in $A((t))^*$ is just the reduction of
Kato's residue (see [12] and [13]) of some lifts of $f$ and $g$. The
existence and the basic properties of the latter follows at once from the
properties of the rigid analytic regulator. Hence the aforementioned
results follow immediately from some well-known facts such as the
deformation theory of smooth projective curves is unobstructed and Weil's
reciprocity law. One may even give a new proof of the Weil reciprocity
law using the observed continuity of the tame symbol by degenerating the
curve to a stable curve with rational components in its special fiber.
Since the Anderson-Romo reciprocity law implies them, the reciprocity
laws of Tate and Witt also follow from Weil's law.
\enddefinition
\definition{Contents 1.9} The goal of the next chapter is to review the
construction of the rigid analytic regulator and to list its basic properties for rational subdomains of the projective line without proofs, mainly for the sake of the reader. The relationship with Kato's residue homomorphism is established in the third chapter. The fourth chapter is concerned with the Contou-Carr\`ere symbol and the Anderson-Romo
reciprocity law. We review the general construction of the residue homomorphism for K\"ahler differentials in the fifth chapter which we relate to the Contou-Carr\`ere symbol in case of local Artinian rings, and to Kato's residue homomorphism in case of local fields of dimension two. This reciprocity law is used to deduce Theorem 1.5 and Corollary 1.7 in chapter six.
\enddefinition
\definition{Acknowledgment 1.10} I wish to thank the IH\'ES for
its warm hospitality and the pleasant environment they created for
productive research, where this article was written.
\enddefinition

\heading 2. Review of the rigid analytic regulator
\endheading

\definition{Notation 2.1} In this chapter all claims are stated without
proof. The interested reader is kindly asked to consult [17]. Let $F$ be
a local field and let $\Bbb C$ denote the completion of the algebraic closure of $F$ with respect to the unique extension of the absolute value on $F$. Recall that $\Bbb C$ is an algebraically closed field complete with respect to an ultrametric absolute value which will be denoted by $|\cdot|$. Let $|\Bbb C|$ denote the set of values of the latter. Let $\Bbb
P^1$  denote the projective line over $\Bbb C$. For any $x\in\Bbb P^1$ and
any two rational non-zero functions $f$, $g\in\Bbb C((t))$ on the
projective line let $\{f,g\}_x$ denote the tame symbol of the pair
$(f,g)$ at $x$. Recall that a subset $U$ of $\Bbb P^1$ is a connected
rational subdomain, if it is non-empty and it is the complement of the
union of finitely many pair-wise disjoint open discs. Let
$\partial U$ denote the set of these complementary open discs. The
elements of $\partial U$ are called the boundary components of $U$, by
slight abuse of language. Let $\Cal O(U)$, $\Cal R(U)$ denote the algebra
of holomorphic functions on $U$ and the subalgebra of restrictions of
rational functions, respectively.  Let $\Cal O^*(U)$, $\Cal R^*(U)$ denote
the group of invertible elements of these algebras. The group $\Cal
R^*(U)$ consists of rational functions which do not have poles or zeros
lying in $U$. For each $f\in\Cal O(U)$ let $\|f\|$ denote $\sup_{z\in
U}|f(z)|$. This is a finite number, and makes $\Cal O(U)$ a Banach
algebra over $\Bbb C$. We say that the sequence $f_n\in\Cal O(U)$
converges to $f\in\Cal O(U)$, denoted by $f_n\rightarrow f$, if $f_n$
converges to $f$ with respect to the topology of this Banach algebra, i.e.
$\lim_{n\rightarrow\infty}\|f-f_n\|=0$. For every real number
$0<\epsilon<1$ we define the sets $\Cal O_{\epsilon}(U)=\{f\in\Cal
O(U)|\|1-f\|\leq\epsilon\}$, and $U_{\epsilon}=\{z\in\Bbb C|
|1-z|\leq\epsilon\}$. 
\enddefinition
\proclaim{Theorem 2.2} There is a unique map $\{\cdot,\cdot\}_D:\Cal
O^*(U)\times\Cal O^*(U)\rightarrow\Bbb C^*$ for every $D\in\partial U$,
called the rigid analytic regulator, with the following properties:

\noindent (i) For any two $f$, $g\in\Cal R^*(U)$ their regulator is:
$$\{f,g\}_D=\prod_{x\in D}\{f,g\}_x,$$

\noindent (ii) the regulator $\{\cdot,\cdot\}_D$ is bilinear
in both variables,

\noindent (iii) the regulator $\{\cdot,\cdot\}_D$ is alternating: 
$\{f,g\}_D\cdot\{g,f\}_D=1$,

\noindent (iv) if $f$, $1-f\in\Cal O(U)^*$, then $\{f,1-f\}_D$ is $1$,

\noindent (v) for each $f\in\Cal O_{\epsilon}(U)$ and $g\in\Cal O^*(U)$
we have $\{f,g\}_D\in U_{\epsilon}$.\ $\square$
\endproclaim
\definition{Remark 2.3} It is an immediate consequence of property $(v)$
that the rigid analytic regulator is continuous with respect to the
supremum topology. Explicitly, if $f$ and $g$ are elements of $\Cal
O^*(U)$, $D\in\partial U$ is a boundary component, and $f_n\in\Cal
O^*(U)$, $g_n\in\Cal O^*(U)$ are sequences such that $f_n\rightarrow f$
and $g_n\rightarrow g$, then the limit
$$\lim_{n\rightarrow\infty}\{f_n,g_n\}_D,$$
exists, and it is equal to $\{f,g\}_D$. 
\enddefinition
Let $\Cal M(U)$ denote the field of meromorphic functions of
$U$ and let $\Cal M^*(U)$ denote the multiplicative group of non-zero elements of $\Cal M(U)$.
\proclaim{Theorem 2.4} There is a unique set of homomorphisms $\deg_D:\Cal M^*(U)\rightarrow\Bbb Z$ where $U$ is any connected rational subdomain and $D\in\partial U$ is a boundary component with the following properties:
\roster
\item"$(i)$" the homomorphism $\deg_D$ is zero on $\Cal O_1(U)$,
\item"$(ii)$" for every $f\in\Cal R^*(U)$ the integer $\deg_D(f)$ is the number of zeros $z$ of $f$ with $z\in D$ counted with multiplicities  minus the number of poles $z$ of $f$ with $z\in D$ counted
with multiplicities,
\item"$(iii)$" for every $f\in\Cal M^*(U)$ we have $\deg_D(f|_Y)=\deg_D(f)$ where $Y\subseteq U$ is any connected rational subdomain satisfying the property $D\in\partial Y$.\ $\square$
\endroster
\endproclaim
\definition{Definition 2.5} If $U$ is still a connected rational subdomain
of $\Bbb P^1$, and $f$, $g$ are two meromorphic functions on $U$, then
 for all $x\in U$ the functions $f$ and $g$ have a power series expansion
around $x$, so in particular their tame symbol $\{f,g\}_x$ at $x$ is
defined. The tame symbols extends to a homomorphism
$\{\cdot,\cdot\}_x:K_2(\Cal M(U))\rightarrow\Bbb C^*$. We define the
group $K_2(U)$ as the kernel of the direct sum of tame symbols:
$$\bigoplus_{x\in U}\{\cdot,\cdot\}_x:K_2(\Cal M(U))\rightarrow
\bigoplus_{x\in U}\Bbb C^*.$$
Let $k=\sum_if_i\otimes g_i\in K_2(U)$, where $f_i$, $g_i\in\Cal M(U)$,
and let $D\in\partial U$. Let moreover $Y$ be a connected rational
subdomain of $U$ such that $f_i$, $g_i\in\Cal O^*(Y)$ for all $i$ and
$\partial U\subseteq\partial Y$. Define the rigid analytical regulator
$\{k\}_D$ by the formula:
$$\{k\}_D=\prod_i\{f_i|_Y,g_i|_Y\}_D.$$
\enddefinition
\proclaim{Theorem 2.6} $(i)$ For each $k\in K_2(U)$ the rigid analytical
regulator $\{k\}_D$ is well-defined, and it is a homomorphism
$\{\cdot\}_D:K_2(U)\rightarrow\Bbb C^*$,

\noindent $(ii)$ for any two functions $f$, $g\in\Cal O^*(U)$ we have
$\{f\otimes g\}_D=\{f,g\}_D$,

\noindent $(iii)$ for every $k\in K_2(U)$ the product of all regulators
on the boundary components of $U$ is equal to 1:
$$\prod_{D\in\partial U}\{k\}_D=1\text{.\ $\square$}$$
\endproclaim
\definition{Definition 2.7} For every $U\subset\Bbb P^1$ connected rational subdomain let $\Bbb Z\partial U$ denote the free abelian group with the elements of $\partial U$ as free generators. Let $H_1(U)$ denote the quotient of $\Bbb Z\partial U$ by the $\Bbb Z$-module generated by $\sum_{D\in\partial U}D$. For every $D\in\partial U$ we let $D$ denote the class of $D$ in $H_1(U)$ as well. Let $\Cal Ab$ denote the category of abelian groups. Let $\Cal Crs$ denote the category whose objects are connected rational subdomains of $\Bbb P^1$ and whose morphisms are holomorphic maps between them. Let $D(w,r)$ denote the open disc of radius $r$ and center $w$, that is
$$D(w,r)=\{z\in\Bbb C||z-w|<r\}$$
where $0<r\in|\Bbb C|$. Finally for every pair $a\leq b$ of numbers in $|\Bbb C|$ let $A(a,b)$ denote the closed annulus $\Bbb P^1-D(0,a)-D(\infty,1/b)$. Of course it is a connected rational subdomain.
\enddefinition
\proclaim{Theorem 2.8} There is a unique functor $H_1:\Cal Crs\rightarrow\Cal Ab$ with the following properties:
\roster
\item"$(i)$" for every $U\subset\Bbb P^1$ connected rational subdomain $H_1(U)$ is the group defined in 3.7,
\item"$(ii)$" for every map $U\rightarrow Y$ which is the restriction of a projective linear transformation $f$ and $D\in\partial U$ boundary component we have:
$$H_1(f)(D)=f(D)\in H_1(Y),$$
\item"$(iii)$" for every $f:U\rightarrow A(a,b)$ holomorphic map and $D\in\partial U$ boundary component we have:
$$H_1(f)(D)=\deg_D(f)D(0,a)\in H_1(A(a,b))\text{.\ $\square$}$$
\endroster
\endproclaim
\definition{Definition 2.9} Let $U\subset\Bbb P^1$ be a connected rational subdomain. For every class $c\in H_1(U)$ and element $k\in K_2(U)$ we define the regulator $\{k\}_c$ as
$$\{k\}_c=\prod_{D\in\partial U}\{k\}_D^{c(D)},$$
where $\sum_{D\in\partial U}c(D)D$ is a lift of $c$ in $\Bbb Z\partial U$. By claim $(iii)$ of Theorem 2.6 this regulator is well-defined. For every holomorphic map $h:U\rightarrow Y$ between two connected rational subdomains let $h^*:K_2(\Cal M(Y))\rightarrow K_2(\Cal M(U))$ be the pull-back homomorphism induced by $h$. By restriction it induces a homomorphism $K_2(Y)\rightarrow K_2(U)$.
\enddefinition
\proclaim{Theorem 2.10} For any $k\in K_2(Y)$ and $c\in H_1(U)$ we
have:
$$\{h^*(k)\}_c=\prod_{E\in\partial Y}\{k\}_{H_1(h)(c)}\text{.\ $\square$}$$
\endproclaim
\definition{Definition 2.11} Let $U$ be now a connected rational subdomain
of $\Bbb P^1$ defined over $F$. This means that
$$U=\{z\in\Bbb P^1||f_i(z)|\leq1\ (\forall i=1,\ldots,n)\}$$
as a set for some natural number $n$ and rational functions $f_1,\ldots,f_n\in F(t)$. Let $\Cal O_F(U)$, $\Cal R_F(U)$, $\Cal O_F^*(U)$, $\Cal R_F^*(U)$ and $\Cal M_F(U)$ denote the algebra of holomorphic functions, the subalgebra of restrictions of $F$-rational functions, the groups of invertible elements of these algebras and the field of meromorphic functions on the rigid analytic space $U$, respectively. Let $U'$ denote the underlying rational subdomain over $\Bbb C$. Let $K_2(U)$ denote the largest subgroup of $K_2(\Cal M_F(U))$ which maps into
$K_2(U')$ under the restriction homomorphism $K_2(\Cal M_F(U))\rightarrow
K_2(\Cal M(U'))$. An $F$-rational boundary component of $U$ is a set
$D\in\partial U$ such that $D$ is the image of the open disc of radius
$1$ and center $0$ under an $F$-linear projective linear transformation of
$\Bbb P^1$.
\enddefinition 
\proclaim{Proposition 2.12} Let $D$ be an $F$-rational boundary component
of $U$, and let $k\in K_2(U)$. Then $\{k\}_D\in F^*$.\ $\square$
\endproclaim

\heading 3. Kato's residue homomorphism
\endheading

\definition{Definition 3.1} In this chapter next we will continue to use
the notation of the second chapter. In this chapter we will also assume that the absolute value on $F$ is induced by a discrete valuation. Let $D$ denote the open disc $D(0,1)$. Let $M$ be the field of fractions of $\Cal O[[z]]$ and let $\widehat M$ denote the completion of $M$ with respect to the discrete valuation of $M$ defined by the prime ideal $\goth m\Cal O[[z]]$ of height one, where $\goth m$ is the unique maximal ideal of $\Cal O$. The field $\widehat M$ is just the field of bidirectional formal Laurent series of the form $\sum_{n\in\Bbb
Z}a_nz^n$ over $F$ such that $|a_n|$ is bounded above and
$\lim_{n\rightarrow-\infty}|a_n|=0$. It is a local field equipped with
the absolute value
$$\|\sum_{n\in\Bbb Z}a_nz^n\|_s=\max_{n\in\Bbb Z}|a_n|.$$
Every element of the formal Laurent series ring $\Cal O[[z]]$ defines a
holomorphic function on the rigid analytic space $D$, hence every element
$M$ gives a meromorphic function on $D$. By Weierstrass' preparation
theorem each element of $\Cal O[[z]]$ is the product of a polynomial and
a unit of this ring, hence it has only a finite number of zeros in $D$.
Therefore the limit
$$\{f,g\}_D=\lim\Sb\epsilon\rightarrow0\\0<\epsilon<1\\
\epsilon\in|\Bbb C|\endSb
\{f,g\}_{D(0,1-\epsilon)}=\prod_{x\in D}\{f,g\}_x$$
becomes stationary for any pair of elements $f,g\in M^*$ and defines an
$F^*$-valued bilinear map satisfying the Steinberg relation by Theorem
2.2 and Proposition 2.12. Therefore it induces a homomorphism
$\{\cdot\}_D:K_2(M)\rightarrow F^*$. Note that the rigid analytic
regulator denoted by the same symbol has the same value as this pairing
for those pairs of functions for which both of them are defined by
Theorem 2.6. Hence our notation will not cause confusion.
\enddefinition
\proclaim{Proposition 3.2} There is a unique homomorphism
$\{\cdot\}_D:K_2(\widehat M)\rightarrow F^*$, called Kato's residue
homomorphism, such that
\roster
\item"$(i)$" the composition of the natural homomorphism $K_2(M)
\rightarrow K_2(\widehat M)$ and Kato's residue homomorphism
is the homomorphism $\{\cdot\}_D$ defined above,
\item"$(ii)$" for each $f\in\widehat M^*$ and $g\in\Cal O[[z]]$ with
$\|1-g\|_s<\epsilon<1$ we have $\{f,g\}_D\in U_{\epsilon}$.
\endroster
\endproclaim
\definition{Proof} Clearly Kato's residue homomorphism is unique if it
exits. We claim that for each $f,g \in M^*$ with $\|1-g\|_s<\epsilon<1$
we have $\{f,g\}_D\in U_{\epsilon}$. We first show that this claim
implies the proposition. In this case we may define
$\{f,g\}_D$ for any two elements $f$ and $g$ of $\widehat M^*$ as the
limit
$$\lim_{n\rightarrow\infty}\{f_n,g_n\}_D,$$
where $f_n\in M^*$, $g_n\in M^*$ are sequences such that
$f_n\rightarrow f$ and $g_n\rightarrow g$. This limit exits because the
sequence above is Cauchy by the claim above. Its value is non-zero as
$$1=\lim_{n\rightarrow\infty}\{f_n,g_n\}_D\cdot\{g_n,f_n\}_D=
\lim_{n\rightarrow\infty}\{f_n,g_n\}_D\cdot
\lim_{n\rightarrow\infty}\{g_n,f_n\}_D.$$
It is also independent of the sequences chosen as any two sequences may
be combed together to show that they give the same limit. The map
$\{\cdot,\cdot\}_D$ defined this way is automatically a  bilinear map
satisfying claim $(ii)$ and the Steinberg relation, hence the existence
follows. For every $1>\delta\in|\Bbb C|$ sufficiently close to
$1$ the holomorphic functions $f$ and $g$ are elements of $\Cal
O^*(A_{\delta})$, where $A_{\delta}$ is the annulus $A_{\delta}=
\{z\in\Bbb C||z|=\delta\}$. Write $1-g=\sum_{n\in\Bbb Z}a_nz^n$ as an
element of $\widehat M$. This power series will converge for all $z\in
A_{\delta}$ when $\delta$ sufficiently close to $1$, hence there is a
number $0<\rho<1$ and a negative integer $N$ such that
$|a_n|\leq\epsilon\rho^{-n}$ for all $n<N$. For all $\delta\in|\Bbb C|$
such that $\rho<\delta<1$ we have the following estimate for the supremum norm $\|1-g\|$ on the annulus $A_{\delta}$: 
$$\|1-g\|\leq\max(\|\sum_{n<N}a_nz^n\|,\|\sum_{n\geq N}a_nz^n\|)\leq
\max(\epsilon,\epsilon\delta^N)=\epsilon\delta^N.$$
Therefore the limit inferior of the supremum norms $\|1-g\|$ on the annuli
$A_{\delta}$ is at most $\|1-g\|_s$, so the claim is now
clear by $(v)$ of Theorem 2.2.\ $\square$
\enddefinition
Let $t\in\Cal O[[z]]$ be a uniformizer, which here means that $t$ is of
the form $cz+z^2h$, where $c\in\Cal O^*$ and $h\in\Cal O[[z]]$. Then there
is a unique $\Cal O$-algebra automorphism $\phi:\Cal O[[z]]\rightarrow\Cal
O[[z]]$ such that $\phi(z)=t$ and $\|\phi(h)\|_s=\|h\|_s$ for every $h\in
\Cal O[[z]]$. The automorphism $\phi$ extends uniquely to a
norm-preserving automorphism $\phi:\widehat M\rightarrow\widehat M$.
Let $\phi_*: K_2(\widehat M)\rightarrow K_2(\widehat M)$ denote the
induced automorphism.
\proclaim{Proposition 3.3} The automorphism $\phi_*$ leaves Kato's residue
homomorphism invariant.
\endproclaim
\definition{Proof} By continuity we only have to show that the
equation $\{\phi_*(k)\}_D=\{k\}_D$ holds for any $k\in K_2(M)$. Note that
the power series $t$ as a holomorphic function $t:D\rightarrow D$ leaves
the annulus $A_{\delta}$ invariant for any positive rational
$\delta<1$ where we continue to use the notation of the proof above. In
fact for any $z\in A_{\delta}$ we have
$$|t(z)-cz|\leq|z^2|=\delta^2<\delta=|cz|.$$
The inequality above also implies that $\deg_{D(0,\delta)}(t)=1$ by claim $(iii)$ of Theorem 2.8, hence the claim follows at once from
Theorem 2.10.\ $\square$
\enddefinition
\proclaim{Lemma 3.4}  For every pair of positive integers $n$ and $m$ the
following identities hold:
\roster
\item"$(i)$" $\{1-at^{-n},1-bt^{-m}\}_D=1$, if $|a|<1$ and $|b|<1$,
\item"$(ii)$" $\{1-at^n,1-bt^m\}_D=1$, if $|a|\leq1$ and $|b|\leq1$,
\item"$(iii)$" $\{1-at^n,1-bt^{-m}\}_D=(1-a^{m/(n,m)}
b^{n/(n,m)})^{(n,m)}$, if $|a|\leq1$ and $|b|<1$,
\item"$(iv)$" $\{1-at^{-n},1-bt^m\}_D=(1-a^{m/(n,m)}
b^{n/(n,m)})^{-(n,m)}$, if $|a|<1$ and $|b|\leq1$.
\endroster
\endproclaim
\definition{Proof} Note that the equations $(iii)$ and $(iv)$ are
equivalent, because we can get the latter from the former by reversing
the roles of the symbols $a$ and $b$, and using the antisymmetry of the
rigid regulator. Hence we only have to show $(i)$, $(ii)$ and $(iii)$. First assume that both $m$ and $n$ are equal to
$1$. We may assume that both $a$ and $b$ are non-zero. In case $(i)$ the
two linear expressions in $t^{-1}$ both have one zero, which are $a$ and
$b$, respectively. They also have a pole, which is the point zero. These points all lie in $D$, so Weil's reciprocity law implies:
$$\{1-at^{-1},1-bt^{-1}\}_D=\prod_{x\not\in D}
\{1-at^{-1},1-bt^{-1}\}^{-1}_x=1.$$
In case $(ii)$ the zeros of the two linear polynomials are $1/a$ and
$1/b$, respectively, which do not lie in $D$. Hence the equation holds in
this case. In case $(iii)$ the expression of $1-at$ does not have a zero
or a pole in $D$, but $1-bt^{-1}$ does, hence:
$$\{1-at,1-bt^{-1}\}_D=\{1-at,1-bt^{-1}\}_0\cdot
\{1-at,1-bt^{-1}\}_b=1-ab.$$
Now assume that $n$ and $m$ are relatively prime and none of them is
divisible by the characteristic of $\Bbb C$. Let $\epsilon_1$,
$\epsilon_2$ be a primitive $n$-th and a primitive $m$-th root of unity,
respectively. Let $\alpha$ and $\beta$ be an $n$-th root of $a$ and an
$m$-th root of $b$, respectively. Since $|\epsilon_1^i\alpha|^n=|a|$ and
$|\epsilon_2^j\beta|^m=|b|$, the condition of claim $(iii)$ hold for these
values, so we get:
$$\{1-at^n,1-bt^{-m}\}_D=\prod_{i=1}^n\prod_{j=1}^m
\{1-\epsilon_1^i\alpha t,1-\epsilon_2^j\beta t^{-1}\}_D=\prod_{i=1}^n
\prod_{j=1}^m(1-\epsilon_1^i\alpha\epsilon_2^j\beta)=1-ab.$$
The other two claims follow similarly. Now assume that $n$ and $m$ are
still relatively prime, but one of them, for example $n$, is divisible by
$p$, the characteristic of $\Bbb C$. In this case $1-at^{\pm n}=(1-\alpha
t^{\pm n/p})^p$, where $\alpha^p=a$, so the claims follow from what we
have proved already, by induction on the exponent of $p$ in the primary
factorization of $n$. In the general case we have:
$$\{1-at^{\pm n},1-bt^{\pm m}\}_D=\{1-at^{\pm n/(n,m)},
1-bt^{\pm m/(n,m)}\}^{(n,m)}_D,$$
which follows from applying Theorem 2.10 to the map $t\mapsto
t^{(n,m)}$.\ $\square$
\enddefinition
\proclaim{Lemma 3.5} For every pair of integers $n$ and $m$ the
following identities hold:
\roster
\item"$(i)$" $\{at^n,bt^m\}_D=(-1)^{nm}a^mb^{-n}$, if both $a$ and $b$ are
non-zero,
\item"$(ii)$" $\{at^n,1-bt^m\}_D=1$, if $a\neq0$, $|b|\leq1$ and
$m$ is positive, or $a\neq0$, $|b|<1$ and $m$ is negative.
\endroster 
\endproclaim
\definition{Proof} In case $(i)$ both expressions have at most one
singularity on the disc $D$ which is the point zero. Therefore
$$\{at^n,bt^m\}_D=\{at^n,bt^m\}_0=(-1)^{nm}a^mb^{-n}.$$
In case $(ii)$ we may immediately reduce to the case $b\neq0$ and $|m|=1$
using the same arguments as the proof above. If $m=1$ then linear
expression $1-bt$ has no singularity on the disc $D$, hence
$$\{at^n,1-bt\}_D=\{at^n,1-bt\}_0=1.$$
In the other case the expression $1-bt^{-1}$ has two singularities on the
disc $D$: a pole at $0$ and a zero at $b$. Therefore
$$\{at^n,1-bt^{-1}\}_D=\{at^n,1-bt^{-1}\}_0\{at^n,1-bt^{-1}\}_b
=1\text{.\ $\square$}$$
\enddefinition
\definition{Definition 3.6} Fix a uniformizer $\pi\in F$ and let $\Cal R$ denote the valuation ring of $\widehat M$. For every $u\in\Cal R$ let $\overline
u\in\Cal R/\pi\Cal R$ denote the reduction of $u$ modulo the proper maximal ideal of $\Cal R$. Note that $\Cal R/\pi\Cal R$ is a local field since it is canonically isomorphic to $\bold f((\overline z))$ where $\bold f$ is the residue field of $F$. Let $\nu$ denote the valuation of
$\Cal R/\pi\Cal R$ normalised such that $\nu(\overline z)=1$. Every element $u\in\widehat M^*$ can be written uniquely in the form $\pi^nv$ for some $n\in\Bbb Z$ and $v\in\Cal R^*$. We define $\deg(u)$ as $\nu(\overline v)$.
\enddefinition
\proclaim{Lemma 3.7} We have $\{c,u\}_D=c^{\deg(u)}$ for every $c\in F^*$
and $u\in\widehat M^*$.
\endproclaim
\definition{Proof} By the continuity and the bilinearity of Kato's residue homomorphism we only have to show that the equation in the claim above is true when $u=dv$ where $d\in F^*$ and $v\in\Cal O[[z]]$. Because $\{c,d\}=1$
by definition we may assume that $u=v$. In this case the number of zeros of the convergent power series $u$ on $D$ counted with multiplicities is exactly $\deg(u)$ so the claim holds.\ $\square$
\enddefinition

\heading 4. The Contou-Carr\`ere symbol and the Anderson-Romo reciprocity
law
\endheading

\definition{Notation 4.1} Let $\bold k$ be a perfect
field and let $\Cal C$ denote the category of local Artinian rings
with residue field $\bold k$. By slight abuse of notation we will let
$\Cal C$ denote the class of objects of this category as well.
\enddefinition
\proclaim{Lemma 4.2} Assume that $\bold k$ has characteristic zero. Then
for every object $A$ in $\Cal C$ there is a homomorphism $i:\bold k
\rightarrow A$ such that the composition of the reduction map $A
\rightarrow\bold k$ modulo the maximal ideal of $A$ and $i$ is the identity
map.
\endproclaim
\definition{Proof} This is a special case of Proposition 6 of [21] on pages
33-34.\ $\square$ 
\enddefinition
\proclaim{Lemma 4.3} Assume that $\bold k$ has positive characteristic p. Then for every object $A$ in $\Cal C$ there is a homomorphism $i:\Bbb W(\bold k)\rightarrow A$ of local rings, where $\Bbb W(\bold k)$ is the ring of Witt
vectors of $\bold k$ of infinite length, such that the map induced by $i$
on the residue fields is the identity.
\endproclaim
\definition{Proof} By Theorem 8 of [21] on page 43 the ring $\Bbb W(\bold k)$ is strict hence the claim follows from Proposition 10 of [21] on pages 38-39.\ $\square$
\enddefinition
\proclaim{Proposition 4.4} Assume that $\bold k$ is algebraically closed
and let $\Cal D$ be a subclass of the class of
objects of $\Cal C$ such that the following conditions hold:
\roster
\item"$(i)$" if $A\in\Cal C$ is the quotient of a discrete valuation ring
with residue field $\bold k$, then $A\in\Cal D$,
\item"$(ii)$" if $A\in\Cal C$ is the quotient of an element of $\Cal D$,
then $A\in\Cal D$,
\item"$(iii)$" if $A\in\Cal C$ and for every $x\in A^*$ different from
$1$ there is a $B\in\Cal D$ and a homomorphism $\phi:A\rightarrow B$ such
that $\phi(x)\neq1$, then $A\in\Cal D$.
\endroster
In this case $\Cal D$ is the whole class of objects of $\Cal C$.
\endproclaim
\definition{Proof} For every pair of natural numbers $n$ and $m$ let
$A_{n,m}$ denote the local Artinian algebra:
$$A_{n,m}=\bold k[[x_1,\ldots,x_n]]/(
\{\prod_{j=1}^nx_j^{J(j)}|
J:\{1,2,\ldots,n\}\rightarrow\Bbb N
,\sum_{j=1}^nJ(j)=m+1\}).$$
First assume that $\bold k$ has characteristic zero. In this case for
every $A\in\Cal C$ there is a surjective local homomorphism
$A_{n,m}\rightarrow A$ for some $n$ and $m$. On the other hand
$A_{1,m}\in\Cal D$ by condition $(i)$. Therefore it will be enough to
show that for every $x\in A^*_{n,m}$ with $x\neq1$ there is a homomorphism
$\phi:A_{n,m}\rightarrow A_{1,m}$ such that $\phi(x)\neq1$ by condition
$(iii)$. There is a positive integer $k\leq m$ such that $x\equiv1\mod\goth m^{k-1}$, but $x\not\equiv1\mod\goth m^k$.
Every local homomorphism $\phi:A_{n,m}\rightarrow A_{1,m}$
induces a $\bold k$-linear homomorphism $T^l_{\phi}:\goth m^l/
\goth m^{l+1}\rightarrow\goth n^l/\goth n^{l+1}$ for every positive $l\leq m$, where $\goth n$ is the maximal proper ideal of $A_{1,m}$. For every vector space $V$ over $\bold k$ let $Sym^l(V)$ denote the $l$-th symmetric power of $V$ and for every $\bold k$-linear map $h:V\rightarrow W$ between vector spaces over $\bold k$ let $Sym^l(h):Sym^l(V)\rightarrow Sym^l(W)$ denote the $l$-th symmetric power of this homomorphism. The multiplication induces a natural isomorphism between $Sym^k(\goth m/\goth m^2)$, $Sym^k(\goth n/\goth n^2)$ and $\goth m^k/\goth m^{k+1}$ and $\goth n^k/\goth n^{k+1}$,
respectively, and under these identifications we have
$T^k_{\phi}=Sym^k(T^1_{\phi})$. Since any
$\bold k$-linear map $h:\goth m/\goth m^2\rightarrow\goth n/\goth
n^2\cong\bold k$ is induced by a local homomorphism
$\phi:A_{n,m}\rightarrow A_{1,m}$, it will be sufficient to prove the
following lemma.
\enddefinition
\proclaim{Lemma 4.5} For every $0\neq v\in Sym^k(\bold k^n)$ there is a
$\bold k$-linear map $\phi:\bold k^n\rightarrow\bold k$ such that
$Sym^k(\phi)(v)\neq0$.
\endproclaim
\definition{Proof} We are going to prove the claim by induction on $n$.
The case $n=1$ is obvious. Let $x_1,x_2,\ldots,x_n$ be a basis of
$\bold k^n$. Write $v$ as
$$v=\sum_{j=0}^kp_j(x_1,x_2,\ldots,x_{n-1})x_n^j,$$
where $p_j\in Sym^{k-j}(\bold k^{n-1})$. For a $0\leq j\leq k$ the
polynomial $p_j$ is not zero, therefore there is a $\bold k$-linear map
$\phi_1:\bold k^{n-1}\rightarrow\bold k$ by induction, where $\bold
k^{n-1}$ is spanned by $x_1,x_2,\ldots,x_{n-1}$, such that
$Sym^{k-j}(\phi_1)(p_j)\neq0$. The polynomial:
$$p(t)=\sum_{j=0}^k Sym^{k-j}(\phi_1)(p_j)t^j$$
is not identically zero, hence it has finitely many roots. The field
$\bold k$ is assumed be algebraically closed, in particular it is not
finite. Hence there is a $\beta\in\bold k$ which is not a root of
the polynomial above. Let $\phi:\bold k^n\rightarrow\bold k$ be the
unique $\bold k$-linear extension of $\phi_1$ with $\phi(x_n)=\beta$. In
this case $Sym^k(v)=p(\beta)\neq0$, so the claim is proved. \ $\square$ 
\enddefinition
Now assume that $\bold k$ has characteristic $p>0$ and let $\Bbb W(\bold
k)$ denote the ring of Witt vectors of $\bold k$ of infinite length. For
every pair of natural numbers $n$ and $m$ let $B_{n,m}$ denote the local
Artinian algebra:
$$B_{n,m}=\Bbb W(\bold k)[[x_1,\ldots,x_n]]/(
\{p^{J(0)}\cdot\prod_{j=1}^nx_j^{J(j)}|
J:\{0,1,\ldots,n\}\rightarrow\Bbb N
,\sum_{j=0}^nJ(j)>m\}).$$
For every $A\in\Cal C$ there is a surjective local homomorphism
$B_{n,m}\rightarrow A$ for some $n$ and $m$. By repeating the argument
above we may reduce the proof of the proposition to show that there is a
homomorphism $\phi:B_{n,m}\rightarrow B_{0,m}$ such that $\phi(x)\neq1$
for every $x\in B^*_{n,m}$ with $x\equiv1\mod\goth m^{k-1}$, but
$x\not\equiv1\mod\goth m^k$ for some positive integer $k\leq m$, where
$\goth m$ is the maximal ideal of $B_{n,m}$. Every local homomorphism
$\phi:B_{n,m}\rightarrow B_{0,m}$
induces a $\bold k$-linear homomorphism $T^l_{\phi}:\goth m^l/\goth
m^{l+1}\rightarrow\goth n^l/\goth n^{l+1}$ for any positive $l\leq m$,
where $\goth n$ is the maximal proper ideal of $B_{0,m}$. Let $T_p$,
$T^{\bot}$ denote the $\bold k$-linear subspace of $\goth m/\goth m^2$ generated by $p$ and the elements $x_1,x_2,\ldots,x_n$, respectively. The multiplication
induces a natural isomorphism:
$$\goth m^k/\goth m^{k+1}=\bigoplus_{j=0}^kSym^j(T_p)\otimes
Sym^{k-j}(T^{\bot})$$
and another between $Sym^k(\goth n/\goth n^2)$ and $\goth n^k/\goth
n^{k+1}$. Moreover there is a canonical isomorphism
$\iota:T_p\rightarrow\goth n/\goth n^2$ between these one-dimensional
vector spaces. Under these identifications we have
$$T^k_{\phi}=\bigoplus_{j=0}^kSym^j(\iota)\otimes
Sym^{k-j}(T^1_{\phi}|_{T^{\bot}}).$$
Since every $\bold k$-linear map $h:T^{\bot}\rightarrow\goth
n/\goth n^2\cong\bold k$ is induced by a local homomorphism
$\phi:B_{n,m}\rightarrow B_{0,m}$, the proposition follows from the
lemma above.\ $\square$
\definition{Definition 4.6} Let $A\in\Cal C$ be a local Artinian ring
with maximal ideal $\goth m$ and let $f$ be any element of $A((t))^*$.
Then there is an integer $w(f)\in\Bbb Z$, and a sequence of elements
$a_i\in A$ indexed by the integers such that $a_0\in A^*$, $a_{-i}\in
\goth m$ for $i>0$, $a_{-i}=0$ for $i$ sufficiently large, and
$$f=a_0\cdot t^{w(f)}\cdot\prod_{i=1}^{\infty}(1-a_it^i)\cdot
\prod_{i=1}^{\infty}(1-a_{-i}t^{-i}),$$
and these are uniquely determined by $f$. The integer $w(f)$ is called the
winding number of $f$ and the elements $a_i$ are called the Witt
coordinates of $f$. Let $f$, $g\in A((t))^*$ be arbitrary with winding
numbers $w(f)$, $w(g)$ and Witt coordinates $a_i$, $b_j$, respectively.
 By definition the Contou-Carr\`ere symbol $\langle f,g\rangle$ is:
$$\langle f,g\rangle=(-1)^{w(f)w(g)}
{a_0^{w(g)}\prod_{i=1}^{\infty}\prod_{j=1}^{\infty}
(1-a_i^{j/(i,j)}b_{-j}^{i/(i,j)})^{(i,j)}\over b_0^{w(f)}
\prod_{i=1}^{\infty}\prod_{j=1}^{\infty}(1-a_{-i}^{j/(i,j)}
b_i^{i/(i,j)})^{(i,j)}}\in A^*.$$
Obviously all but finitely many terms are equal to one in the infinite
product above, hence the Contou-Carr\`ere is a well-defined alternating
map:
$$\langle\cdot,\cdot\rangle:A((t))^*\times A((t))^*\rightarrow A^*.$$
It is also clear from the formula that the Contou-Carr\`ere symbol is
equal to the tame symbol if $A$ is a field. 
\enddefinition
\proclaim{Proposition 4.7} The Contou-Carr\`ere symbol is a bilinear map
satisfying the Steinberg relation.
\endproclaim
\definition{Proof} Because $\bold k$ is perfect for every object $A$ of $\Cal
C$ there is a local Artinian ring $B$ with residue field $\overline{\bold k}$, where the latter is the algebraic closure of $\bold k$, and an injective
local homomorphism $i:A\rightarrow B$. Indeed the algebra $B=A\otimes_{\bold
k}\overline{\bold k}$ and $B=A\otimes_{\Bbb W(\bold k)}\Bbb W(\overline
{\bold k})$ will do, when $\bold k$ has characteristic zero or positive characteristic,
respectively, using the fact has $A$ can be equipped with the structure of
a $\bold k$-algebra or $\Bbb W(\bold k)$-algebra, respectively, by Lemmas
4.2 and 4.3. Hence we may assume that $\bold k$
is algebraically closed. Let $\Cal D$ denote the subclass of those local
Artinian rings with residue field $\bold k$ which satisfy the claim
of the proposition above. Clearly we only have to show that this subclass
satisfies the conditions of Proposition 4.4. If $A\in\Cal C$ is the
quotient of a discrete valuation ring $R$ with residue field $\bold k$,
we may assume that $R$ is complete with respect to its valuation. Let $K$
be the quotient field of $R$ and let $\Bbb C$ be the completion of the
algebraic closure of $K$. The latter is an algebraically closed field
complete with respect to an ultrametric absolute value. For every $f\in
A((t))^*$ there is a lift $\tilde f\in R((t))^*$ whose image is $f$ under
the functorial map $R((t))\rightarrow A((t))$. By Lemmas 3.4 and 3.5 the
Contou-Carr\`ere symbol of $f$ and $1-f$ is the reduction of rigid
analytic regulator $\{\tilde f,1-\tilde f\}_D\in R$ modulo the maximal
ideal of $R$, hence the Contou-Carr\`ere symbol satisfies the Steinberg
relation. A similar argument shows that it is also bilinear, therefore
$(i)$ of Proposition 4.4 holds for $\Cal D$. Property $(ii)$ also follows
from same reasoning, because every $f\in A((t))^*$ has a lift $\tilde
f\in B((t))^*$ if the map $B\rightarrow A$ is surjective. Finally let
$A\in\Cal C$ be an algebra which satisfies the condition in $(iii)$ of
Proposition 4.4. Assume that there is an $1\neq f\in A((t))^*$ such that
$\langle f,1-f\rangle\neq1$. Then there is a $B\in\Cal D$ and a
homomorphism $\phi:A\rightarrow B$ such that $1\neq\phi(\langle
f,1-f\rangle) =\langle\phi_*(f),1-\phi_*(f)\rangle\neq1$, where
$\phi_*:A((t))\rightarrow B((t))$ is the functorial map induced by $\phi$,
which is a contradiction. A similar argument shows that the
Contou-Carr\`ere symbol is bilinear over $A$, therefore property $(iii)$
also holds for $\Cal D$.\ $\square$
\enddefinition
Let $x\in A[[t]]$ be an uniformizer, which means that $x$ is of the form
$ct+t^2h$, where $c\in A^*$ and $h\in A[[t]]$. In this case there is a
unique $A$-algebra automorphism $\phi:A[[t]]\rightarrow A[[t]]$ such that
$\phi(t)=x$. On the other hand every $A$-algebra automorphism of
$A[[t]]$ is of this form. The automorphism $\phi$ extends uniquely to an
automorphism $\phi:A((t))\rightarrow A((t))$ by localizing at the maximal
ideal.
\proclaim{Proposition 4.8} The automorphism $\phi$ leaves the
Contou-Carr\`ere symbol in\-vari\-ant.
\endproclaim
\definition{Proof} As in the proof above we may assume that $\bold k$
is algebraically closed. Let $\Cal D$ again denote the subclass of those local Artinian rings with residue field $\bold k$ which satisfy the claim
of the proposition above. We need to show only that this subclass
satisfies the conditions of Proposition 4.4. Let $\psi:B\rightarrow A$ be
a surjective homomorphism of local algebras with residue field $\bold k$
and let $\psi_*:B((t))\rightarrow A((t))$ be the functorial map induced
by $\psi$. If $B$ is Artinian or a discrete valuation ring then there is a $B$-algebra automorphism $\phi_B:B((t))\rightarrow B((t))$ such that
$\psi_*\circ\phi_B=\phi\circ\psi_*$ which is of type described before
Proposition 3.3 if $B$ is a discrete valuation ring. Hence $(i)$ and
$(ii)$ of Proposition 4.4 holds for $\Cal D$ by Proposition 3.3. A similar
argument as above shows that condition $(iii)$ also holds for $\Cal D$.\ $\square$
\enddefinition
\definition{Notation 4.9} Let $A\in\Cal C$ be a local Artinian ring
and let $\pi:X\rightarrow\text{Spec}(A)$ be a projective flat morphism
whose fiber $X_0$ over the unique closed point of Spec$(A)$ is a reduced,
connected, regular curve over $\bold k$. Let $S$ be a finite set of
closed points of $X$ (or, equivalently, of $X_0$) and let $f$ and $g$ be
two rational functions on $X$ which are invertible on the complement of
$S$. For every $s\in S$ choose an $A$-algebra isomorphism $\phi_s$
between the completion $\widehat{\Cal O}_{s,X}$ of the stalk ${\Cal O}_{s,X}$
of the structure sheaf of
$X$ at $s$ and $A[[t]]$. The latter induces an isomorphism between the
localization of $\widehat{\Cal O}_{s,X}$ by the semigroup of those
elements whose image under the canonical map $\widehat{\Cal O}_{s,X}
\rightarrow\widehat{\Cal O}_{s,X_0}$ is non-zero, where $\widehat{\Cal O}_{s,X_0}$ is the  the completion of the stalk ${\Cal O}_{s,X_0}$, and $A((t))$, which will be denoted by $\phi_s$, by the usual abuse of notation. Let $\langle
f,g\rangle_s$ denote the Contou-Carr\`ere symbol of the image of $f$ and
$g$ under $\phi_s$ for every $s$ in $S$. By Proposition 4.8 the value of
$\langle f,g\rangle_s$ is independent of the choice of the isomorphism
$\phi_s$, so the symbol $\langle f,g\rangle_s$ is well-defined. The
following result is the reciprocity law of Anderson and Romo (see [2]).
\enddefinition
\proclaim{Proposition 4.10} The product of all the Contou-Carr\`ere
symbols of $f$ and $g$ is equal to 1:
$$\prod_{s\in S}\langle f,g\rangle_s=1.$$
\endproclaim
\definition{Proof} We are going to use the same strategy for proof as
we used before: in particular we assume that $\bold k$
is algebraically closed. Let $\Cal D$ denote the subclass of those local Artinian rings with residue field $\bold k$ which satisfy the claim of the proposition above. We will show that this subclass satisfies the
conditions of Proposition 4.4. If $A\in\Cal C$ is the quotient of a
discrete valuation ring $R$ with residue field $\bold k$, we may again
assume that $R$ is complete with respect to the valuation. Let $K$ denote
the quotient field of $R$ and let $\Bbb C$ denote the completion of the
algebraic closure of $K$ as above. Because the deformation theory of
regular projective curves is unobstructed, there is a formal scheme
$\goth X$ over the formal spectrum of $R$ whose fiber over Spec$(A)$ is
$X$. By the algebraicity theorem of Grothendieck $\goth X$ is actually the
formal completion of a smooth curve over Spec$(R)$ which will be denoted
by $\goth X$ by abuse of notation. By flatness there are rational
functions $\tilde f$ and $\tilde g$ on $\goth X$ whose restriction to
the  fiber over Spec$(A)$ is $f$ and $g$, respectively. The rigid
analytic domain $D_s$ of $\Bbb C$-valued points of $\goth X$ which reduces
to $s$ is isomorphic to the open disc $D$ by the formal inverse function
theorem.  By Lemmas 3.4 and 3.5 the Contou-Carr\`ere symbol of $f$ and
$g$ is the reduction of rigid analytic regulator the product of the tame
symbols $\{\tilde f,\tilde g\}_x$ modulo the maximal ideal of $R$ where
$x$ is running through the $\Bbb C$-valued points of set $D_s$. The
rational functions $\tilde f$ and $\tilde g$ has only poles or zeros in the
union of the sets $D_s$ hence the reciprocity law of Anderson and
Romo holds by Weil's reciprocity law.

Property $(ii)$ also follows from same reasoning. If the map $B\rightarrow
A$ is surjective and $X$, $f$ and $g$ are as above, then there is a
similar triple $\tilde X$, $\tilde f$ and $\tilde g$ over Spec$(B)$ such
that $X$ is the fiber of $\tilde X$ over Spec$(A)$ is $X$ and the
restriction of $\tilde f$ and $\tilde g$ to $X$ is $f$ and $g$,
respectively, because the deformation theory of $X$ is unobstructed. If
$B\in\Cal D$ then the claim holds for the triple $\tilde X$, $\tilde f$,
$\tilde g$, so it must hold for the triple $X$, $f$, $g$ as well. Finally
let $A\in\Cal C$ be an algebra which satisfies the condition in $(iii)$
of Proposition 4.4. Assume that there are rational functions $f$ and $g$
as above such that $\prod_{s\in S}\langle f,g\rangle_s\neq1$. Then there
is a $B\in\Cal D$ and a homomorphism $\phi:A\rightarrow B$ such that
$1\neq\phi(\prod_{s\in S}\langle
f,g\rangle_s)=\prod_{s\in\phi^*(S)}\langle\phi^*(f),
\phi^*(g)\rangle_s$, where $\phi^*(f)$, $\phi^*(g)$ and $\phi^*(S)$ are
the base change of the corresponding objects on the curve $\phi^*(X)$
which is the base change of $X$ with respect to the map
$\phi^*:\text{Spec}(B)\rightarrow\text{Spec}(A)$. This is a contradiction,
therefore property $(iii)$ also holds for $\Cal D$.\ $\square$
\enddefinition
\definition{Remark 4.11} It is possible to push the methods of this paper a bit further to actually give a proof of Weil's reciprocity law itself
by reducing it to the case of Mumford curves, when it follows from
$(iii)$ of Theorem 2.6 at once. We will only sketch this argument because
it uses a considerable amount of machinery compared to the relatively
elementary nature of Weil's reciprocity law. For any scheme $S$ and any
stable curve $\pi:C\rightarrow S$ of genus $g$ let $\omega_{C/S}$ denote
relative dualizing sheaf. By Theorem 1.2 of [8], page 77 the functor
which assigns to each scheme $S$ the set of stable curves
$\pi:C\rightarrow S$, and an isomorphism $\Bbb P(\pi_*(\omega_{C/S}^
{\otimes3}))\cong\Bbb P^{5g-6}_S$ (modulo isomorphism) is represented by a
fine moduli scheme $\goth H_g$. By Corollary 1.7 of [8], page 83 and the
main result of [8], pages 92-96, the scheme $\goth H_g$ is smooth over
the spectrum of $\Bbb Z$ and the base change $(\goth H_g)_{\text{\rm
Spec}(\bold k)}$ is irreducible for any algebraically closed field $\bold
k$. Let $X$ smooth, projective curve over $\bold k$ and let $f$ and $g$ be
two non-zero rational functions on $X$. We may assume that the genus $g$
of $X$ is at least two by taking a cover of $X$ and proving the
reciprocity law for the pull-back of $f$ and $g$ instead. Let $x$ be a
$\bold k$-valued point of $\goth H_g$ such that the underlying curve is
$X$ and let $y$ be another $\bold k$-valued point such that the underlying
curve is totally degenerate. Since $(\goth H_g)_{\text{\rm Spec}
(\bold k)}$ is an irreducible, smooth quasi-projective variety, 
repeated application of Bertini's theorem shows that there is a smooth,
irreducible curve $S$ mapping to $(\goth H_g)_{\text{\rm Spec}(\bold k)}$ whose image contains both $x$ and $y$. Let $\pi:C\rightarrow S$ be the pull-back of the universal family. There are rational functions $\tilde f$ and $\tilde g$ on $C$ whose restriction to the fiber over $x$, which is
$X$, are $f$ and $g$, respectively. Since the base change of $C$ to the spectrum of the local field of $S$ at $y$ is a Mumford curve, Weil's reciprocity law holds for $\tilde f$ and $\tilde g$, hence holds for $f$
and $g$, too. One may say that this proof is close in spirit to the
classical proof of the reciprocity law over the complex numbers using
triangulation, since it decomposes the curve to small pieces in a suitable
topology.
\enddefinition

\heading 5. The differential reciprocity law
\endheading

\definition{Definition 5.1} We will continue to use the notation of the previous chapter. For every $\bold k$-algebra $A$ let $\Omega_A^{\cdot}$ denote the graded differential algebra of $\bold k$-linear K\"ahler differential forms of $A$ and for every $\bold k$-algebra homomorphism $h:A\rightarrow B$ let $\Omega^k(h):\Omega^{\cdot}_A\rightarrow\Omega^{\cdot}_B$ induced by $h$ by functoriality. Every $\omega\in\Omega^k_{A((t))}$ can be written uniquely in the form:
$$\omega=\sum_{i=1}^m\beta_i{dt\over t^i}+\omega_0$$
where $m$ is a natural number, $\beta_i\in\Omega^{k-1}_A$ and $\omega_0\in\Omega^k_{A[[t]]}+A((t))\Omega^k_A$. Let $\text{Res}^k(\omega)\in\Omega^{k-1}_A$ denote the element $\beta_1$. We get a map
$\text{Res}^k:\Omega^k_{A((t))}\rightarrow\Omega^{k-1}_A$ which is called the residue.
\enddefinition
\proclaim{Proposition 5.2} The following holds:
\roster
\item"$(i)$" we have $\text{\rm Res}^{k+i}(\alpha\omega)=\alpha\text{\rm Res}^k
(\omega)$ for every $\alpha\in\Omega^i_A$ and $\omega\in\Omega^k_{A((t))}$,
\item"$(ii)$" we have $\Omega^{k-1}(h)\circ\text{\rm Res}^k=\text{\rm Res}^k\circ
\Omega^k(h')$ where $h:A\rightarrow B$ is a $\bold k$-algebra homomorphism and $h':A((t))\rightarrow B((t))$ is the corresponding $\bold k$-algebra homomorphism induced by functoriality,
\item"$(iii)$" we have $\text{\rm Res}^k(\omega)=0$ for every $\omega\in\Omega^k_{A[[t]]}$ and for every $\omega\in\Omega^k_{A[{1\over t}]}$,
\item"$(iv)$" the map $\text{\rm Res}^k$ does not depend on choice of the uniformizer $t$.
\endroster
\endproclaim
\definition{Proof} Our method of proving the first two claims is the same. Using the notation of Definition 5.1 we have:
$$\alpha\omega=\sum_{i=1}^m\alpha\beta_i{dt\over t^i}+\alpha\omega_0.$$
Because $\alpha\beta_i\in\Omega^{k+i-1}_A$ and $\alpha\omega_0\in\Omega^{k+i}_{A[[t]]}+A((t))\Omega^{k+i}_A$ we have $\text{\rm Res}
^{k+i}(\alpha\omega)=\alpha\beta_1$ by definition so claim $(i)$ is true. On the other hand:
$$\Omega^k(h')(\omega)=\sum_{i=1}^m\Omega^{k-1}(h)(\beta_i)
{dt\over t^i}+\Omega^k(h')(\omega_0)$$
where $\Omega^{k-1}(h)(\beta_i)\in\Omega^{k-1}_B$ and $\Omega^k(h')(\omega_0)\in\Omega^k_{B[[t]]}+B((t))\Omega^k_B$. Therefore we have $\text{\rm Res}^k(\Omega^k(h')(\omega))=\Omega^{k-1}(h)(\beta_1)$ as claim $(ii)$ says. The first half of claim $(iii)$ is immediate from the definition of the residue. In order to prove the second half we only need to show the identity $\text{Res}^1(dt^{-n})=0$ for all $n\geq1$ by the $\Omega_A^{\cdot}$-linearity of the residue spelled out in claim $(i)$. But the latter is obvious. Claim $(iv)$ means the following: let $x\in A[[t]]$ be an uniformizer, which means that $x$ is of the form $tu$, where $u\in A[[t]]^*$. In this case there is a unique $A$-algebra automorphism $\phi:A[[t]]\rightarrow A[[t]]$ such that $\phi(t)=x$ as we already saw when we prepared to formulate Proposition 4.8. The automorphism $\phi$ extends uniquely to an
automorphism $\phi:A((t))\rightarrow A((t))$ by localizing at the maximal
ideal. Claim $(iv)$ means that the equation $\text{\rm Res}^k\circ\Omega^k(\phi)
=\text{\rm Res}^k$ holds. Because the homomorphism $H^k(\phi)$ maps $\Omega^k_{A[[t]]}$ and $A((t))\Omega^k_A$ into itself we only need to show that
$$\text{Res}^1({dx\over x})=1\text{\ and \ }\text{Res}^1({dx\over x^n})=0
\text{\ for all $n\geq2$}$$
by $\Omega_A^{\cdot}$-linearity. These identities follow at once from the same type of identity in Proposition 5' in [20] on pages 20-21 by the principle of prolongation of algebraic identities quoted in the proof of the proposition just mentioned above. (Or one may use the functorial argument explained in the remark following the proof of Proposition 5' in [20] instead.)\ $\square$
\enddefinition
\definition{Notation 5.3} For every $\bold k$-algebra $B$ let $dlog:B^*\rightarrow
\Omega^1_B$ denote the logarithmic differential given by the rule $dlog(u)=u^{-1}du$ for every $u\in B^*$. Let moreover $dlog^2$ denote the $\Bbb Z$-bilinear pairing:
$$dlog^2:B^*\otimes B^*\rightarrow\Omega^2_B$$
given by the rule:
$$dlog^2(a,b)=dlog(a)dlog(b)\quad(\forall a\in B^*,\forall b\in B^*).$$
Recall that for every object $A$ of $\Cal C$ the symbol $\langle\cdot,\cdot\rangle$ denotes the Contou-Carr\`ere symbol.
\enddefinition
\proclaim{Proposition 5.4} The following diagram commutes:
$$\CD A((t))^*\otimes A((t))^*@>dlog^2>>\Omega^2_{A((t))}\\
@V\langle\cdot,\cdot\rangle VV@VV\text{\rm Res}^2V\\A^*@>dlog>>
\Omega^1_A\endCD$$
for every object $A$ of $\Cal C$.
\endproclaim
\definition{Proof} By bilinearity and antisymmetry of the Contou-Carr\`ere symbol and the map $dlog^2$ it will be sufficient to prove for every pair of integers $n$, $m$ and element $a$, $b\in A$ the following identities:
\roster
\item"$(i)$" $\text{Res}^2(dlog^2(1-at^n,1-bt^m))=0$, if $n$, $m>0$,
\item"$(ii)$" $\text{Res}^2(dlog^2(1-at^{-n},1-bt^{-m}))=0$, if $a$, $b\in\goth m$ and $n$, $m>0$,
\item"$(iii)$" $\text{Res}^2(dlog^2(at^n,1-bt^m))=0$, if $a\in A^*$ and $m>0$,
\item"$(iv)$" $\text{Res}^2(dlog^2(at^n,1-bt^m))=0$, if $a\in A^*$, $b\in\goth m$ and $m<0$,
\item"$(v)$" $\text{Res}^2(dlog^2(at^n,bt^m))=m(da/a)-n(db/b)$, if $a$, $b\in A^*$,
\item"$(vi)$" $\text{Res}^2(dlog^2(1-at^n,1-bt^{-m}))=dlog((1-a^{m/(n,m)}
b^{n/(n,m)})^{(n,m)})$, if $b\in\goth m$ and $n$, $m>0$,
\item"$(vii)$" $\text{Res}^2(dlog^2(1-ft^{Mn+1},1-bt^{-n}))=0$, if $f\in A[[t]]$, $b\in\goth m$ and $n>0$,
\endroster
where $\goth m$ is the maximal proper ideal of $A$ and $M$
is a positive integer such that $\goth m^M=0$. Note that
$$dlog(1-at^n)=-(dat^n+nat^{n-1}dt)(1+at^n+a^2t^{2n}+\cdots)$$
lies in $\Omega^1_{A[[t]]}$, if $a\in A$ and $n>0$, and lies in
$\Omega^1_{A[{1\over t}]}$, if $a\in\goth m$ and $n<0$. Hence the first two identities follow from claim $(iii)$ of Proposition 5.2. For every $a\in A^*$ and $b\in A$ we have:
$$\split dlog^2(at^n,1-bt^m)=&-(1+bt^m+b^2t^{2m}+\cdots)({da\over a}+n{dt\over t})
(dbt^m+mbt^{m-1}dt)\\
=&(ndb-{mb\over a}da)(t^{m-1}dt+bt^{2m-1}dt+b^2t^{3m-1}dt+\cdots)+\omega_0\endsplit$$
where $\omega_0\in A((t))\Omega^2_A$ when either $m>0$ or when $b\in\goth m$ and $m<0$.
The first summand in the second line lies in $\Omega^2_{A[[t]]}$, if $m>0$, and lies in $\Omega^2_{A[{1\over t}]}$, if $b\in\goth m$ and $m<0$.
Hence its residue is zero so identities $(iii)$ and $(iv)$ are true. For every $a$, $b\in A^*$ we have:
$$dlog^2(at^n,bt^m)=({da\over a}+n{dt\over t})({db\over b}+m{dt\over t})=
(m{da\over a}-n{db\over b}){dt\over t}+\omega_0$$
where $\omega_0\in\Omega^2_A$ so identity $(v)$ is clear. By definition:
$$dlog((1-a^{m/(n,m)}b^{n/(n,m)})^{(n,m)})=-
{ma^{m/(n,m)}b^{n/(n,m)}{da\over a}+na^{m/(n,m)}b^{n/(n,m)}{db\over b}
\over1-a^{m/(n,m)}b^{n/(n,m)}}$$
for every $a\in A$, $b\in\goth m$ and $n$, $m>0$. We also have:
$$(1-at^n)^{-1}(1-b^{-m})^{-1}=\sum_{k\in\Bbb Z}\sum\Sb i,j\in\Bbb N\\in-jm=k\endSb
a^ib^jt^k$$
for such $a$ and $b$. Because $in-jm=m-n$ for any $i,j\in\Bbb N$ if and only if $i+1=
lm/(n,m)$ and $j+1=ln(n,m)$ for some $l\in\Bbb N$ we have:
$$(1-at^n)^{-1}(1-b^{-m})^{-1}=(ab)^{-1}\big(\sum_{l=1}^{\infty}a^{lm/(n,m)}b^{ln/(n,m)}
\big)t^{m-n}+(r+s)t^{m-n}$$
for some $r\in{1\over t}A[{1\over t}]$ and $s\in tA[[t]]$.
Hence we have:
$$\split dlog^2(1-at^n,1-bt^{-m})=&
{(dat^n+nat^{n-1}dt)(dbt^{-m}-mbt^{-m-1}dt)\over(1-at^n)(1-bt^{-m})}\\
=&{-(mbda+nadb)t^{n-m-1}dt\over(1-at^n)(1-bt^{-m})}+\omega_0\\
=&-{a^{m/(n,m)-1}b^{n/(n,m)-1}(mbda+nadb)t^{-1}dt\over1-a^{m/(n,m)}b^{n/(n,m)}}
+\omega_0+\omega_1\endsplit$$
where $\omega_0\in A((t))\Omega^2_A$ and $\omega_1\in\Omega^2_{A[{1\over t}]}
+\Omega^2_{A[[t]]}$. Identity $(vi)$ is now obvious. Finally consider the
last identity. Note that
$$dlog(1-bt^{-n})=-(dbt^{-n}-nbt^{-n-1}dt)(1+bt^{-n}+b^2t^{-2n}+\cdots+b^{M-1}
t^{(1-M)n})$$
because $b^M=0$ by assumption, and
$$dlog(1-ft^{Mn+1})=-(dft^{Mn+1}+(Mn+1)ft^{Mn}dt)(1+ft^{Mn+1}+f^2t^{2Mn+2}+\cdots)$$
hence
$$\split dlog^2(1-ft^{Mn+1},1-bt^{-n})=&
(dft+(Mn+1)fdt)(db-nbt^{-1}dt)g\\
=&(tdfdb-nbdfdt+(Mn+1)fdtdb)g\endsplit$$
where $g\in A[[t]]$. The claim is now clear.\ $\square$
\enddefinition
\definition{Definition 5.5} Let $L$ be a field complete with respect to
a discrete valuation and let $\Cal R$, $\goth m$ denote its discrete valuation ring and the maximal ideal of $\Cal R$, respectively. Assume that the residue
field of $L$ is $\bold k$ and the quotient map $\Cal R\rightarrow\bold k$
has a section which is a ring homomorphism. The latter equips $L$ and $\Cal
R$ with a $\bold k$-algebra structure. Let $\widehat{\Omega}^{\cdot}_L$ denote the graded differential algebra which is the quotient of the complex $\Omega^{\cdot}_L$ by the homogeneous ideal generated by $\cap_{n\geq1}\goth m^n\Omega^{\cdot}_{\Cal R}$ and let $\widehat{\Omega}^k_{\Cal R}$ denote the image of $\Omega^k_{\Cal R}$ in $\widehat{\Omega}^d_L$ under the quotient map. For every natural number $n$ let $\Cal R_n$ denote the truncated ring $\Cal R/\goth m^{n+1}$
and for every pair $m\leq n$ of natural numbers let $\pi_n:\Cal R
\rightarrow\Cal R_n$ and $\pi_{n,m}:\Cal R_n\rightarrow\Cal R_m$ denote the canonical projection. The system modules $\{\Omega^k_{\Cal R_n}\}_{n\in\Bbb N}$ form a compatible system with respect to the morphisms $\Omega^k(\pi_{n,m})$ ($m\leq n$) hence it has a projective limit $\varprojlim_{n\rightarrow\infty}
(\Omega^k_{\Cal R_n})$. The maps
$\Omega^k(\pi_n):\Omega^k_{\Cal R}\rightarrow\Omega^k_{\Cal R_n}$ factor through $\widehat{\Omega}^k_{\Cal R}$ and their limit induces an identification: $\widehat{\Omega}^k_{\Cal R}
\cong\varprojlim_{n\rightarrow\infty}(\Omega^k_{\Cal R_n})$ which we will use without further notice. Let $dlog:L^*\rightarrow
\widehat{\Omega}^1_L$ and  $dlog^2:K_2(L)\rightarrow\widehat{\Omega}^2_L$
also denote the composition of $dlog$, $dlog^2$ an the quotient map
$\Omega^1_L\rightarrow\widehat{\Omega}^1_L$, $\Omega^2_L\rightarrow
\widehat{\Omega}^1_L$, respectively. 
\enddefinition
\definition{Definition 5.6} Let $\pi$ be a uniformizer of $L$ and let
$\widehat{\Omega}^k_L(log)$ denote the subgroup $\pi^{-1}\widehat{\Omega}^k_{\Cal
R}$ of $\widehat{\Omega}^k_L$. Clearly the group $\widehat{\Omega}^k_L(log)$
is independent of the choice of the uniformizer $\pi$. Let $\Cal O$ denote the discrete valuation ring $\bold k[[x]]$ and let $F$ denote its quotient field. Let $\widehat M$ denote the field attached to $F$ which was introduced in Definition 3.1 and let $\Cal R$ denote the valuation ring of $\widehat M$. The uniformizer $x$ of $F$ is also a unifomizer in $\widehat M$. There is a natural isomorphism $\Cal R_n\cong\Cal O_n((\overline z))$ for every
$n\in\Bbb N$, where $\overline z$ denotes also the reduction of $z$ in $\Cal R_n$ for every $n$ by slightly extending the notation introduced in Definition
3.6, therefore for every $\omega\in\Omega^k_{\Cal R_n}$ the
residue $\text{Res}^k(\omega)\in\Omega^{k-1}_{\Cal O_n}$ is well-defined. For every $\omega\in\widehat{\Omega}^k_{\widehat M}(log)$ let $\text{Res}^k
(\omega)\in\widehat{\Omega}^{k-1}_F$ be given by the rule:
$$\text{Res}^k(\omega)={1\over x}
\varprojlim_{n\rightarrow\infty}
(\text{Res}^k(\widehat{\Omega}^k(\pi_n)(x\omega)))$$
where the map $\widehat{\Omega}^k(\pi_n):\widehat{\Omega}^k_{\Cal R}
\rightarrow\Omega^k_{\Cal R_n}$ is induced by $\Omega^k(\pi_n)$. The system:
$$\{\text{Res}^k(\widehat{\Omega}^k(\pi_n)(x\omega))\}_{n\in\Bbb N}$$
satisfies the compatibility described above by claim $(ii)$ of Proposition 5.2 hence $\text{Res}^k(\omega)$ is well-defined. Because of the $\Cal O_n$-linearity of the residue it is obvious that $\text{Res}^k(\omega)$ is independent of the choice of $x$ as the notation indicates.
\enddefinition
\definition{Remark 5.7} Let $\phi:\widehat M\rightarrow\widehat M$ be a valuation-preserving
$F$-algebra automorphism. Then there is a unique map $\widehat{\Omega}^k(\phi):
\widehat{\Omega}^k_{\widehat M}
\rightarrow\widehat{\Omega}^k_{\widehat M}$ such that $\widehat{\Omega}^k(\phi)
\circ q_k=q_k\circ\Omega^k(\phi)$ where $q_k:\Omega^k_{\widehat M}
\rightarrow\widehat{\Omega}^k_{\widehat M}$ is the quotient map. The automorphism
$\widehat{\Omega}^k(\phi)$ of $\widehat{\Omega}^k_{\widehat M}$ preserves the subgroup $\widehat{\Omega}^k_{\widehat M}(log)$ and it commutes with the residue map $\text{Res}^k$ by claim $(iv)$ of Proposition 5.2.
\enddefinition
\proclaim{Theorem 5.8} We have $dlog^2(k)\in \widehat{\Omega}^2_{\widehat
M}$ for every $k\in K_2(\widehat M)$ and the diagram:
$$\CD K_2(\widehat M)@>dlog^2>>\widehat{\Omega}^2_{\widehat M}(log)\\
@V\{\cdot,\cdot\}_DVV@VV\text{\rm Res}^2V\\F^*@>dlog>>
\widehat{\Omega}^1_F\endCD$$
is commutative where $\{\cdot,\cdot\}_D$ denotes Kato's residue homomorphism.
\endproclaim
\definition{Proof} By the linearity of the $dlog^2$ map we only have to verify
the first claim of the theorem as well as the identity expressed by the commutative
diagram above for the elements of any set of generators of $K_2(\widehat M)$. Hence we may assume that $k=u\otimes v$ where either $u$, $v\in\Cal R^*$ or $u=x$ and $v$ is arbitrary. In the first case we have $dlog^2(k)\in\widehat{\Omega}^k_{\Cal R}$ obviously and the identity holds by Proposition 5.4. In the second case
we may write $v$ in the form $v=x^nw$ for some $n\in\Bbb Z$ and $w\in\Cal R^*$. Because $\{x,x\}_D=1$ and $dlog^2(x\otimes x)=0$ by definition we may assume that $v=w$. The first claim is now obvious. Moreover in this case we may write $v$ in the form $v=z^{\deg(v)}t$ for some $t\in\Cal R^*$ such that the reduction $t_k$ of $t$ modulo $x^k\Cal R$ lies in $\Cal
O_k[[\overline z]]^*\subset\Cal R_k$ for every $k\in\Bbb N$. Therefore $dlog(t_k)\in
\Omega^1_{\Cal O_k[[\overline z]]}$ and we have:
$$\split\text{Res}^2(\Omega^2(\pi_k)(dx{dv\over v}))=&
\text{Res}^2(\deg(v)\Omega^2(\pi_k)(dx{dz\over z}))+
\text{Res}^2(\Omega^2(\pi_k)(dx){dt_n\over t_n})\\
=&\deg(v)\Omega^1(\pi_k)(dx)\endsplit$$
for every $k\in\Bbb N$. The claim now follows from Lemma 3.7.\ $\square$
\enddefinition

\heading 6. The image and kernel of the rigid analytic regulator in positive
characteristic
\endheading

\definition{Notation 6.1} For every scheme $X$ let $K_{2,X}$ denote the sheaf on $X$ associated to the presheaf $U\mapsto K_2(H^0(U,\Cal O_X))$ for the Zariski-topology where $K_2(A)$ denotes Milnor's $K$-group of any ring $A$. Let $\Bbb W_n\Omega^*_X$ denote the de Rham-Witt pro-complex of any ringed topos $X$ of $\Bbb F_p$-algebras. Moreover we let $F$ denote the Frobenius morphism of the de Rham-Witt pro-complex.  Recall that the logarithmic differential $dlog^1:\Cal O_X^*\rightarrow\Bbb W_n\Omega^1_X$ is defined as the composition of the Teichm\"uller lift $\Cal O_X^*\rightarrow\Bbb W_n\Omega^0_X$ and the differential $d:\Bbb W_n\Omega^0_X\rightarrow\Bbb W_n\Omega^1_X$, where $X$ is the same as above. The bilinear map of sheaves:
$$dlog^2:\Cal O_X^*\times\Cal O_X^*@>>>\Bbb W_n\Omega^2_X$$
given by the formula:
$$dlog^2(f\otimes g)=dlog^1(f)dlog^1(g)$$
also satisfies the Steinberg relation $dlog^2(f\otimes(1-f))=0$ for all
$f\in\Cal O_X^*$ with $1-f\in\Cal O_X^*$, hence it induces a map $dlog^2:K_{2,X}
\rightarrow\Bbb W_n\Omega^2_X$. Moreover let $\nu_n(k)$ denote the kernel of $1-F$ on the degree $k$ term $\Bbb W_n\Omega^k_X$ of the de Rham-Witt pro-complex on the topos $X$. Let $\Bbb W_n\Omega^i_{X,log}$ denote the abelian sub-sheaf generated by the image
of $dlog^i$, where $i=1$, $2$. It is easy to see using the defining
relations of the de Rham-Witt pro-complex that $\Bbb
W_n\Omega^i_{X,\log}$ lies in $\nu_n(i)$.
\enddefinition
We will need the following result which is a special case of the celebrated
theorem in [6] due to Bloch, Gabber and Kato.
\proclaim{Theorem 6.2} Let $F$ be a field of characteristic $p$. Then the
map
$$K_2(F)/p^nK_2(F)@>dlog^2>>H^0(F_{et},\Bbb W_n\Omega^2_{F_{et},log})$$
is an isomorphism, where $H^0(F_{et},\Bbb W_n\Omega^2_{F_{et},log})$ denotes the group of global sections of the sheaf $\Bbb W_n\Omega^2_{F_{et},log}$ on the
\'etale site of the spectrum of $F$.
\endproclaim
\definition{Proof} The map is well-defined as $\Bbb W_n\Omega^*$ is
annihilated by $p^n$. The map is an isomorphism by Corollary 2.8 of [6],
page 117-118.\ $\square$
\enddefinition
\definition{Notation 6.3} Let $\bold k$ be a perfect field as in the previous
two chapters. For every $\bold k$-scheme $X$ let $\Omega_X^{\cdot}$ denote the complex of graded differential $\Cal O_X$-algebras of $\bold k$-linear K\"ahler differential forms on $X$. Note that the complex $\Omega_X^{\cdot}$ is canonically isomorphic to the complex $\Bbb W_1\Omega^{\cdot}_X$. In particular there is a map $dlog^2:K_{2,X}\rightarrow\Omega_X^{\cdot}$. For every $k\in\Bbb N$ and for every Cartier divisor $D$ on $X$ let $\Omega^k_X(D)$ denote the sheaf $\Omega^k_X\otimes_{\Cal O_X}\Cal O_X(D)$. Let $i:X-D\rightarrow X$ denote the open immersion of the complement of the support of $D$ into $X$. Then the pull-back $i^*\Omega^k_X(D)$ is canonically isomorphic to $\Omega^k_{X-D}$. Let
$$i^*:H^0(X,\Omega^k_X(D))\rightarrow H^0(X-D,\Omega^k_{X-D})$$
denote the composition of the pull-back and this identification, too.
\enddefinition
\proclaim{Lemma 6.4} Assume that $X$ is a smooth surface over $\bold k$ and $D$ is a normal crossings divisor. Then the image of the map:
$$dlog^2:H^0(X-D,K_{2,X})\rightarrow H^0(X-D,\Omega_{X-D}^2)$$
lies in the image of the map $i^*:H^0(X,\Omega^k_X(D))
\rightarrow H^0(X-D,\Omega^k_{X-D})$ introduced above.
\endproclaim
\definition{Proof} The claim is clearly local on $X$ with respect to the Zariski topology hence we may assume that $X$ is the spectrum of an integral regular $\bold k$-algebra $A$. We may also assume that $D$ has at most one singular point and its branches are the zeros of some elements of $A$. The localization sequence for $K$-theory induces a complex:
$$\CD H^0(X,K_{2,X})@>>>H^0(X-D,K_{2,X})
@>T>>\endCD$$
$$\CD\bigoplus_{C\in\Cal V(\Gamma(D))}H^0(X-S(D),\Cal O^*_C)@>>>\bigoplus_{e\in S(D)}\Bbb Ze\endCD$$
which is exact at the term $H^0(X-D,K_{2,S})$, where the second map is the direct sum of tame symbols along the irreducible components and the third map is the sum of the maps which assigns every element of $H^0(C-S(D),\Cal O^*_C)$ to its divisor considered as a zero cycle supported on $S(D)$ for every $C\in\Cal V(\Gamma(D))$. Let $k$ be an arbitrary element of $H^0(X-D,K_{2,X})$.
Assume first that $D$ is irreducible and let $t\in A$ be an element
whose zero scheme is $D$. Pick an element $u\in A$ whose pull-back to $D$ is equal to the tame symbol $T(k)$. By shrinking $X$ further we may assume that $u\in A^*$. Then $T(k)=T(t\otimes u)$ hence $dlog^2(k-t\otimes u)$ is the pull-back of a differential form on $X$ by the localization sequence. On the other hand $dlog^2(t\otimes u)$ clearly lies in the image of the map
$i^*$. Assume now that $D$ has one ordinary double point $s$ and let $t_1$, $t_2\in A$ be two elements whose zeros are the two branches of $D$. According to the complex above there is an $n\in\Bbb Z$ such that the valuation of the restriction of $T(k)$ onto the zero scheme of $t_1$ and $t_2$ at $s$
is $n$ and $-n$, respectively. Hence by shrinking $X$ further we may assume that there are $u_1$, $u_2\in A^*$ such that the restriction of $T(k)$ onto the zero scheme of $t_1$ and $t_2$ is the restriction of $u_1t_2^n$ and $u_2t_1^{-n}$, respectively. Then we have:
$$T(k)=T(t_1\otimes t_2)^n\cdot T(t_1\otimes u_1)\cdot T(t_2\otimes u_2)$$
and we may argue as above to conclude the proof.\ $\square$
\enddefinition
The lemma above has the following important corollary: because $X-D$ is Zariski-dense
in $X$ the map $i^*:H^0(X,\Omega^k_X(D))
\rightarrow H^0(X-D,\Omega^k_{X-D})$ is injective. Hence the map $dlog^2$
has a unique lift:
$$dlog^2:H^0(X-D,K_{2,X})\rightarrow H^0(X,\Omega_X^2(D))$$
which will be denoted by the same symbol by the usual abuse of notation.
\proclaim{Proposition 6.5} Assume that $\goth X$ is a smooth irreducible projective surface over $\bold k$ and the field $\bold k$ is finite of characteristic
$p$. Then the group $H^0(\goth X,K_{2,\goth X})$ is the extension of a torsion group by its maximal $p$-divisible subgroup.
\endproclaim
\definition{Proof} Using the notation of [15] on pages 307 and 309 let $H^2(\goth X,\Bbb Z(2))$ denote the projective limit $\varprojlim(H^0(\goth X,\nu_n(2)))$. The logarithmic differentials
$$dlog^2:K_2(\goth X)\rightarrow H^0(\goth X,\nu_n(2))$$
satisfy the obvious compatibility hence they induce a map
$$dlog^2:K_2(\goth X)\rightarrow H^2(\goth X,\Bbb Z(2)).$$
Let $\Cal F(\goth U)$ denote the function field of $\goth X$.  Let $\Cal P$ denote the set of prime divisors of $\goth X$ and for every $P\in\Cal P$ let $\bold f_P$ denote the function field of the irreducible curve $P$.
The localization sequence for $K$-theory furnishes an exact sequence:
$$\CD 0@>>>H^0(\goth X,K_{2,\goth X})@>>>K_2(\Cal F(\goth X))
@>T>>\bigoplus_{P\in\Cal P}\bold f_P^*\endCD$$
where the second map is the direct sum of tame symbols along the irreducible components. Every element $k\in H^0(\goth X,K_{2,\goth X})$ of the kernel of $dlog^2$ lies in
$M(\Cal F(\goth X))=\cap_{n\in\Bbb N}p^nK_2(\Cal F(\goth X))$ by the Bloch-Kato-Gabber theorem. Since $K_2(\Cal F(\goth X))$ has no $p$-torsion by Theorem 1.10
of [23] on page 10, the group $M(\Cal F(\goth X))$ is the maximal $p$-divisible subgroup of $K_2(\Cal F(\goth X))$. If the element $l$ is in $M(\Cal F(\goth X))\cap H^0(\goth X, K_{2,\goth X})$ and $k\in M(\Cal F(\goth X))$ is its
unique $p^n$-th root then $T(k)$ is $p^n$-torsion by the localization sequence. But the groups $\bold f_P^*$ has no non-zero $p$-torsion so $k$ lies in the image of $H^0(\goth X,K_{2,\goth X})$. Therefore we get that $M(\Cal F(\goth X))\cap H^0(\goth X, K_{2,\goth X})$, the kernel of the map $dlog^2$ in $H^0(\goth X,K_{2,\goth X})$, is $p$-divisible. Hence it will be sufficient to show that the group $H^2(\goth X,\Bbb Z(2))$ is torsion. This is proved in [15] (the claim itself can be found on page 335) although the proof is somewhat dispersed over the article. It is an immediate consequence of Proposition 5.4 of the paper cited above on page 330-331, the validity of Weil's conjectures for crystalline cohomology (Remark 5.5 of [15] on page 331), and the exact sequence on page 335 of the same
paper.\ $\square$
\enddefinition
\definition{Notation 6.6} Let $A\in\Cal C$ be a local Artinian $\bold k$-algebra and let $\pi:\Bbb P^1_A\rightarrow\text{Spec}(A)$ be the projective line
over $A$. Let $S$ be a finite set of sections $s:\text{Spec}(A)\rightarrow \Bbb P^1_A$ and for every $s\in S$ let $s_0$ denote the $\bold k$-valued point $s_0:\text{Spec}(\bold k)\rightarrow\Bbb P^1_{\bold k}$ we get from
$s$ via base change. Assume that $s_0$ is different from $t_0$ for every
pair $s$, $t\in S$ of different sections. For every $s\in S$ choose an $A$-algebra isomorphism $\phi_s$ between the completion $\widehat{\Cal O}_{s_0,\Bbb P^1_A}$ of the stalk ${\Cal O}_{s_0,\Bbb P^1_A}$ of the structure sheaf of $\Bbb P^1_A$ at $s_0$ and $A[[t]]$. The latter induces an isomorphism between the localization ${\Cal L}_s$ of $\widehat{\Cal O}_{s_0\Bbb P^1_A}$ by the semigroup of those elements whose image under the canonical map
$\widehat{\Cal O}_{s_0,\Bbb P^1_A}\rightarrow\widehat{\Cal O}_{s_0,\Bbb P^1_{\bold
k}}$ is non-zero, where $\widehat{\Cal O}_{s_0,\Bbb P^1_{\bold k}}$ is the completion of the stalk ${\Cal O}_{s_0,\Bbb P^1_{\bold k}}$, and $A((t))$, which will be denoted by $\phi_s$ as well. The image of $s$ is a locally principal closed subscheme of codimension one in $\Bbb P^1_A$ for every element $s$ of $S$. Let $S$ also denote Cartier divisor which is the sum of these divisors by slight abuse of notation. For every $s\in S$ let $\text{Res}^k_s$ denote the composition of the map:
$$H^k_s:H^0(\Bbb P^1_A-S,\Omega^k_{\Bbb P^1_A-S})@>>>\Omega^k_{{\Cal L}_s}
@>\Omega^k(\phi_s)>>\Omega^k_{A((t))},$$
where the first arrow is induced by the tautological map $\Bbb P^1_A-S\rightarrow
\text{Spec}({\Cal L}_s)$, and the residue $\text{Res}^k:\Omega^k_{A((t))}
\rightarrow\Omega^{k-1}_A$. By claim $(iv)$ of Proposition 5.2 the map:
$$\text{Res}^k_s:H^0(\Bbb P^1_A-S,\Omega^k_{\Bbb P^1_A-S})@>>>\Omega^{k-1}_A$$
is independent of the choice of $\phi_s$. Recall that there
is a canonical inclusion $H^0(\Bbb P^1_A,\Omega^2_{\Bbb P^1_A}(S))
\subset H^0(\Bbb P^1_A-S,\Omega^k_{\Bbb P^1_A-S})$.
\enddefinition
\proclaim{Proposition 6.7} The sequence:
$$\CD0\rightarrow\Omega^2_A@>\pi^*>>
H^0(\Bbb P^1_A,\Omega^2_{\Bbb P^1_A}(S))@>\oplus_{s\in S}
\text{\rm Res}^2_s>>\bigoplus_{s\in S}\Omega^1_A
@>\Sigma_{s\in S}(\cdot)>>\Omega^1_A\rightarrow0\endCD$$
is exact where $\pi^*$ is the pull-back with respect to the map
$\pi:\Bbb P^1_A\rightarrow\text{\rm Spec}(A)$.
\endproclaim
\definition{Proof} By base change we may assume that $\bold k$ is algebraically
closed, which implies that it is infinite. Let $R\supseteq S$ be any finite set. Note that for every $\omega\in{1\over t}\Omega^2_{A[[t]]}$ we have $\omega\in\Omega^2_{A[[t]]}$ if and only if $\text{Res}^2(\omega)=0$. Therefore
$$H^0(\Bbb P^1_A,\Omega^2_{\Bbb P^1_A}(S))=
\{\omega\in H^0(\Bbb P^1_A,\Omega^2_{\Bbb P^1_A}(R))|\text{Res}^2_s(\omega)
=0\ (\forall s\in R-S)\}.$$
Hence it is sufficient to prove the proposition above for $R$ instead of $S$. In particular we may assume that the point at infinity $\infty\in\Bbb P^1_A$ lies in $S$ after a suitable automorphism of the $A$-scheme $\Bbb
P^1_A$. Let $x$ be the coordinate function of the affine line $\Bbb A^1_A=\Bbb P^1_A-\infty$. For every $\infty\neq s\in S$ let the same letter denote the
unique element of $A$ such that the image of the section $s$ is the zero
scheme of $x-s\in A[x]$. Every $\omega\in H^0(\Bbb P^1_A-S,\Omega^2_{\Bbb P^1_A-S})$ can be written uniquely in the form:
$$\omega=\omega_0+\sum_{s\in S-\infty}\sum_{k=1}^{n(s)}{\omega_{s,k}\over(x-s)^k}dx+
\sum_{j=0}^{n(\infty)}\omega_{\infty,j}x^jdx$$
where $\omega_0\in\Omega^2_A$, $n(s)$, $n(\infty)\in\Bbb N$ and $\omega_{s,k}$, $\omega_{\infty,j}\in\Omega^1_A$. For every $\infty\neq s\in S$ we may assume
that $x-s$ maps to $t$ with respect to $\phi_s$. Then it is obvious that
$$H^2_s(\eta x^ndx),H^2_s(\eta(x-r)^{-n}dx)\in\Omega^2_{A[[t]]}$$
for every $\eta\in\Omega^1_A$, $n\in\Bbb N$ and $s\neq r\in S-\infty$. Therefore we have $\omega_{s,k}=0$ for every $k>1$ when
$\omega\in H^0(\Bbb P^1_A,\Omega^2_{\Bbb P^1_A}(S))$. We may assume also
that $x^{-1}$ maps to $t$ with respect to $\phi_{\infty}$. In this case it is obvious that
$$H^2_{\infty}(\omega_{s,1}{dx\over x-s})=-\omega_{s,1}{dt\over t}+\eta_s$$
for some $\eta_s\in\Omega^2_{A[[t]]}$ for every $s\in S-\infty$ but
$$H^2_{\infty}(\omega_{\infty,j}x^jdx)=-\omega_{\infty,j}t^{-j-2}dt$$
for every $j=0,1,\ldots,n(\infty)$ so we must have:
$$\omega=\omega_0+\sum_{s\in S-\infty}{\omega_{s,1}\over x-s}dx.$$
By the above $\text{Res}^2_s(\omega)=\omega_{s,1}$ for every $s\in S-\infty$ and $\text{Res}^2_{\infty}(\omega)=-\sum_{s\in S-\infty}\omega_{s,1}$ so the claim is now obvious.\ $\square$
\enddefinition
\definition{Notation 6.8} No we are going to consider the same situation
that we looked at in the introduction. Let $B$ be a smooth irreducible projective curve over $\bold k$ and let $\pi:\goth X\rightarrow B$ be a regular irreducible projective surface fibred over $B$ such that the fiber $\goth X_{\infty}$ of $\goth X$ over the closed point $\infty$ of $B$ is totally degenerate. Then the base change $X$ of $\goth X$ to the completion $F$ of the function field of $B$ with respect to the valuation corresponding to $\infty$ is a Mumford curve over $F$. Let $\goth U\subset\goth X$ be an open subvariety such that its complement is a normal crossings divisor $D$ which is the preimage of a finite set of closed points of $B$ containing $\infty$. The base change of $\goth X$ to the valuation ring of $F$ is a semi-stable model of $X$ whose fiber is $\goth X_{\infty}$ hence the rigid analytic regulator $\{\cdot\}$ introduced in Definition 5.12 of [17] supplies a diagram:
$$\CD H^0(\goth U,,K_{2,\goth U})@>>> H^2_{\Cal M}(X,\Bbb Z(2))
@>\{\cdot\}>>\Cal H(\Gamma(\goth X_{\infty}),F^*),\endCD$$
where the first homomorphism is induced by functoriality. This composition will be denoted by the symbol $\{\cdot\}$ as well. 
\enddefinition
\definition{Definition 6.9} For every $\omega\in H^0(\goth X,
\Omega^2_{\goth X}(D))$ we are going to define a function $\text{Res}(\omega):
\Cal E(\Gamma(\goth X_{\infty}))\rightarrow \widehat{\Omega}^1_F$
as follows. Fix an edge $e\in\Cal E(\Gamma(\goth X_{\infty}))$ and let $s\in S(\goth X_{\infty})$ denote the image of $e$ under the normalization map. Let $C$ be the irreducible component of $\goth
X_{\infty}$ which corresponds to the original vertex of $e$ under the identification
of Notation 1.2. Let $\widehat{\Cal O}_{s,\goth X}$ denote the completion of the stalk $\Cal O_{s,\goth X}$ of the structure sheaf of $\goth X$ at $s$ and let $t\in\Cal O_{s,\goth X}$ be an element whose zero scheme is the germ of the curve $C$. Because $t$ generates a prime ideal in $\widehat{\Cal O}_{s,\goth X}$ the latter gives rise to a discrete valuation on the quotient field $M_e$ of $\widehat{\Cal O}_{s,\goth X}$. Let $\widehat M_e$ denote the completion of $M_e$ with respect to this valuation
and let $i_e:\text{Spec}(\widehat M_e)\rightarrow\goth X$ denote the tautological
map. Note that the closure of the image of the stalk $\Cal O_{\infty,B}$
of the structure sheaf of $B$ at $\infty$ in $\Cal O_{s,\goth X}$ with respect to the map induced by $\pi:\goth X\rightarrow B$ in $\widehat{\Cal O}_{s,\goth
X}$ is isomorphic canonically to the valuation ring $\Cal O$ of $F$. Hence $\widehat M_e$ is canonically equipped with the structure of an $F$-algebra. Let $\phi:\widehat M_e\rightarrow\widehat M$ be the unique valuation-preserving
$F$-algebra homomorphism such that $\phi(t)=x$ where we continue to use
the notation of the previous chapter. Note that $q_2(\Omega(\phi)(i^*_e(\omega)))\in
\widehat{\Omega}^k_{\widehat M}(log)$ where $q_k:\Omega^k_{\widehat M}
\rightarrow\widehat{\Omega}^k_{\widehat M}$ is the quotient map. Hence the value:
$$\text{Res}(\omega)(e)=\text{Res}_2(q_2(\Omega(\phi)(i^*_e(\omega))))\in\widehat{\Omega}^1_F$$
is well-defined and it is independent of the
choice of the element $t$ by Remark 5.7.
\enddefinition
For every oriented graph $G$ and commutative group $R$ let $\Cal F(G,R)$ denote the group of functions $f:\Cal E(G)\rightarrow R$.
\proclaim{Theorem 6.10} We have $\text{\rm Res}(dlog^2(k))\!\in\Cal H(\Gamma(\goth X_0),\widehat{\Omega}^1_F)$ for every $k\!\in H^0(\goth U,K_{2,\goth U})$ and the diagram:
$$\CD H^0(\goth U,K_{2,\goth U})@>\{\cdot\}>>
\Cal H(\Gamma(\goth X_{\infty}),F^*)\\
@V dlog^2VV@VV dlog V\\
H^0(\goth X,\Omega^2_{\goth X}(D))@>\text{\rm Res}>>
\Cal F(\Gamma(\goth X_{\infty}),\widehat{\Omega}^1_F)\endCD$$
is commutative.
\endproclaim
\definition{Proof} We are going to show that $dlog(\{k\}(e))=
\text{Res}(dlog^2(k))(e)$ for every edge $e\in\Cal E(\Gamma(\goth X_{\infty}))$. Then the theorem will follow immediately because $\{k\}$ is a harmonic cochain.
The identity above follows immediately from Theorem 5.8 and the following
alternate description of the rigid analytic regulator. The pull-back of $k$
with respect to $i_e$ is an element $i^*_e(k)\in K_2(\widehat M_e)$. Let
$\phi_*:K_2(\widehat M_e)\rightarrow K_2(\widehat M)$ be the homomorphism
induced by $\phi$. Then we have $\{k\}(e)=
\{\phi_*(i^*_e(k))\}_D$.\ $\square$
\enddefinition
\proclaim{Proposition 6.11} Let $k$ be an element of $H^0(\goth U,
K_{2,\goth U})$ such that $\text{\rm Res}(dlog^2(k))=0$. Then
$dlog^2(k)=0$, too.
\endproclaim
\definition{Proof} Let $x\in F$ be a uniformizer. The closed subscheme of $B$ defined by the $n$-th power of the defining sheaf of ideals of the closed subscheme $\infty$ is isomorphic canonically to Spec$(\Cal O_n)$ where $\Cal O_n=\Cal O/x^n\Cal O$ as in chapter 5. Let $i_n:\text{Spec}(\Cal O_n)\rightarrow B$ be the closed immersion corresponding to this isomorphism. For every irreducible component $C\in\Cal V(\Gamma(\goth X_{\infty}))$ let $C_n$ denote the closed subscheme of $\goth X$ defined by the $n$-th power of the defining sheaf of ideals of the closed subscheme $C$. Let $c_n:C_n\rightarrow\goth X$ be the closed immersion. Then there is a unique morphism $p_n:C_n\rightarrow\text{Spec}(\Cal O_n)$ such that $c_n\circ\pi=p_n\circ i_n$. As an $\Cal O_n$-scheme $C_n$ is isomorphic to the projective line over Spec$(\Cal O_n)$. Let $S$ denote
the Cartier divisor on $C_n$ which is the pull-back of the divisor on $\goth
X$ that is the sum of those irreducible components of $\goth X_{\infty}$
which are intersecting $C$ with respect to the map $c_n$ and are different
from $C$. Then $S$ is the sum of images of sections of the map $p_n$. Let $C_0$ be the divisor of the element $x\in\Cal O_n\subset H^0(C_n,\Cal O_{C_n})$. Multiplication by $x$ induces a map $\Cal O(S+C_0)\rightarrow\Cal O(S)$. By our assumptions the residues of $x c_n^*(dlog^2(k))\in H^0(C_n,\Omega^2_{C_n}(S))$ introduced in Definition 6.6 are all zero. Hence $x c_n^*(dlog^2(k))\in\Omega^2_{\Cal
O_n}$ by Proposition 6.7. But $\Omega^2_{\Cal O_n}=0$ hence we get that the formal completion of $dlog^2(k)$ along the closed scheme $\goth X_{\infty}$
must be zero. The claim is now clear.\ $\square$
\enddefinition
Assume now that $\bold k$ is a field of characteristic $p$.
\proclaim{Corollary 6.12} The map:
$$\{\cdot\}:H^0(\goth U,K_{2,\goth U})/p^nH^0(\goth U,K_{2,\goth U})\rightarrow
H^0(\Gamma(\goth X_{\infty}),F^*/(F^*)^{p^n})$$
induced by the regulator $\{\cdot\}$ is injective for every natural number $n$.
\endproclaim
\definition{Proof} We are going to prove the claim by induction on $n$. Let $\Cal F(\goth U)$ denote the function field of $\goth U$. Assume first that $n=1$ and let $k\in H^0(\goth U,K_{2,\goth U})$ be an element such that $\{k\}\in\Cal H(\Gamma(\goth X_{\infty}),(F^*)^p)$. By Theorem 6.10 we have:
$$\text{Res}(dlog^2(k))=dlog\circ\{k\}=0,$$
hence $dlog^2(k)=0$ by Proposition 6.11. Therefore $k=pl$ for some $l\in K_2(\Cal F(\goth U))$ by the Bloch-Gabber-Kato theorem. Let $\Cal P$ denote the set of prime divisors of $\goth U$ and for every $P\in\Cal P$ let $\bold f_P$ denote the function field of the irreducible curve $P$. The localization sequence for $K$-theory furnishes an exact sequence:
$$\CD 0@>>>H^0(\goth U,K_{2,\goth U})@>>>K_2(\Cal F(\goth U))@>>>\bigoplus_{P\in\Cal P}\bold f_P^*\endCD$$
where the second map is the direct sum of tame symbols along the irreducible components. The image of $l$ with respect to the second map is $p$-torsion. But the groups $\bold f_P^*$ has no non-zero $p$-torsion so $l$ is the image
of an element of $H^0(\goth U,K_{2,\goth U})$. Assume now that claim is proved for $n-1$ and let $k\in H^0(\goth U,K_{2,\goth U})$ be an element such that $\{k\}\in\Cal H(\Gamma(\goth X_{\infty}),(F^*)^{p^n})$. By the induction hypothesis there is an element $l\in H^0(\goth U,K_{2,\goth U})$ such that $k=p^{n-1}l$. Because the group $F^*$ has no $p$-torsion we have $\{l\}\in\Cal H(\Gamma(\goth X_{\infty}),(F^*)^p)$ therefore $l\in pH^0(\goth U,K_{2,\goth U})$ by the above. Hence $k\in p^nH^0(\goth U,K_{2,\goth U})$ as we wished to prove.\ $\square$
\enddefinition
Assume now that $\bold k$ is a finite field of characteristic $p$.
\proclaim{Theorem 6.13} The following holds:
\roster
\item"$(i)$" the quotient group $H^0(\goth U,K_{2,\goth U})/H^0(\goth X,K_{2,\goth X})$ is a finitely generated abelian group whose rank is at most as large as the rank of the group $\Cal H(\Gamma(D),\Bbb Z)$.
\item"$(ii)$" the kernel $\text{\rm Ker}(\{\cdot\})$ of the regulator $\{\cdot\}:
H^0(\goth U,K_{2,\goth U})
\rightarrow\Cal H(\Gamma(\goth X_{\infty}),F^*)$ has a subgroup of finite index which lies in $H^0(\goth X,K_{2,\goth X})$,
\item"$(iii)$" the kernel $\text{\rm Ker}(\{\cdot\})$ above is $p$-divisible.
It is  torsion if Parshin's conjecture holds, and it is finite if the Bass conjecture holds,
\item"$(iv)$" the image $\text{\rm Im}(\{\cdot\})$ of the regulator $\{\cdot\}:H^0(\goth U,K_{2,\goth U})\rightarrow\Cal H(\Gamma(\goth X_{\infty}),F^*)$ is
$p$-sa\-tu\-ra\-ted,
\item"$(v)$" the rank of $\text{\rm Im}(\{\cdot\})$ is at most as large as the rank of the group $\Cal H(\Gamma(D),\Bbb Z)$, 
\item"$(vi)$" the image $\text{\rm Im}(\{\cdot\})$ is discrete if $D=\goth X_{\infty}$.
\endroster
\endproclaim
\definition{Proof} Let us recall that $S(D)$, $\Cal V(\Gamma(D))$ denote the set of singular points and the set of irreducible components of the curve $D$, respectively. The localization sequence for $K$-theory induces a complex:
$$\CD H^0(\goth X,K_{2,\goth X})@>>>H^0(\goth U,K_{2,\goth U})
@>T>>\bigoplus_{C\in\Cal V(\Gamma(D))}
H^0(C-S(D),\Cal O^*_C)\endCD$$
$$\CD@>>>\bigoplus_{e\in S(D)}\Bbb Ze\endCD$$
which is exact at the term $H^0(\goth U,K_{2,\goth U})$, where the second map is the direct sum of tame symbols along the irreducible components and the third map is the sum of the maps which for every $C\in\Cal V(\Gamma(D))$ assigns every element of $H^0(C-S(D),\Cal O^*_C)$ to its divisor considered as a zero cycle supported on $S(D)$. The kernel of the latter is a finitely generated abelian group of rank $\Cal H(\Gamma(D),\Bbb Z)$ hence claim $(i)$ is clear. By Corollary 6.9 the kernel Ker$(\{\cdot\})$ of the map:
$$\{\cdot\}:H^0(\goth U,K_{2,\goth U})\rightarrow\Cal H(\Gamma(\goth X_{\infty}),F^*)$$
is $p$-divisible. Therefore its image with respect to the map $T$ above is finite because the maximal $p$-divisible subgroup of a finitely generated abelian group is finite. Hence the kernel of $T$ in Ker$(\{\cdot\})$ is a subgroup of finite index which lie in $H^0(\goth X,K_{2,\goth X})$. Therefore claim $(ii)$ holds. 

According to Parshin's conjecture the group $H^0(\goth X,K_{2,\goth X})$ should be torsion. Then the same is true for Ker$(\{\cdot\})\cap
H^0(\goth X,K_{2,\goth X})$ and therefore Ker$(\{\cdot\})$ is torsion as well by claim $(ii)$. The Bass conjecture states that $H^0(\goth X,K_{2,\goth X})$ should be a finitely generated abelian group. Hence the same is true for its subgroup Ker$(\{\cdot\})\cap H^0(\goth X,K_{2,\goth X})$. Note that this group is also $p$-divisible: every element of Ker$(\{\cdot\})\cap H^0(\goth X,K_{2,\goth X})$ has a $p$-th root in Ker$(\{\cdot\})\subseteq H^0(\goth U,K_{2,\goth U})$. On the other hand if $px\in H^0(\goth X,
K_{2,\goth X})$ for some $x\in H^0(\goth U,K_{2,\goth U})$ then
$x\in H^0(\goth X,K_{2,\goth X})$ using the localization sequence the same way we did in the proof of Corollary 6.12 already. Hence Ker$(\{\cdot\})\cap H^0(\goth X,K_{2,\goth X})$ is a finite group whose order is not divisible by $p$ so Ker$(\{\cdot\})$ is finite as well by claim $(ii)$. Claim $(iii)$ is now proved.

Because the maximal $p$-divisible subgroup of $\Cal H(\Gamma(\goth X_{\infty}),F^*)$ is finite the image of $H^0(\goth X,K_{2,\goth X})$ with respect to the rigid analytic regulator is torsion by Proposition 6.5. But the torsion of $\Cal H(\Gamma(\goth X_{\infty}),F^*)$ is finite so $\text{Im}(\{\cdot\})$ is finitely
generated and claim $(v)$ is true by claim $(i)$. On the other hand
note that a finitely generated subgroup $\Lambda\subset\Cal H(\Gamma(\goth X_{\infty}),F^*)$ is $p$-saturated if and only if
$$p^n\Lambda=\Lambda\cap p^n\Cal H(\Gamma(\goth X_{\infty}),F^*)$$
for every $n\in\Bbb N$. The latter holds for $\text{Im}(\{\cdot\})$ by Corollary 6.12 so claim $(iv)$ is true. Let
$\text{Reg}:H^0(\goth U,K_{2,\goth U})
\rightarrow\Cal H(\Gamma(\goth X_{\infty}),\Bbb Z)$
denote the tame regulator which is defined as follows. For every $k\in
H^0(\goth U,K_{2,\goth U})$ and for every edge
$e\in\Cal E(\Gamma(\goth X_{\infty}))$ we define Reg$(k)(e)$ as the
valuation of the tame symbol of $k$ along the irreducible component
$o(e)$ of $\goth X_{\infty}$ with respect to the valuation corresponding to the closed point which is the image of $e$ with respect to the normalization map. By Theorem 5.6 of [17] the diagram:
$$\CD H^0(\goth U,K_{2,\goth U})@>\{\cdot\}>>
\Cal H(\Gamma(\goth X_{\infty}),F^*)\\
@VV\text{\rm Reg}V@VVvV\\
\Cal H(\Gamma(\goth X_{\infty}),\Bbb Z)@=
\Cal H(\Gamma(\goth X_{\infty}),\Bbb Z)\\
\endCD$$
commutes where $v$ is the map induced by the normalized valuation on $F$. If $D=\goth X_{\infty}$ then the kernel of Reg contains $H^0(\goth
X,K_{2,\goth X})$ as a subgroup of finite index according to the complex
we wrote down above. Since $H^0(\goth X,K_{2,\goth X})$ is $p$-divisible its image with respect to the regulator $\{\cdot\}$ is finite. Hence the kernel of the map $v$ in $\text{Im}(\{\cdot\})$ is finite, too. Therefore $\text{Im}(\{\cdot\})$ must be discrete as claim $(vi)$ says.\ $\square$
\enddefinition

\enddocument